\documentclass[12pt,a4paper]{article}

\usepackage{amsmath,amsfonts,amsxtra,amssymb,latexsym, amscd,amsthm}
\usepackage[mathscr]{eucal}
\usepackage{indentfirst}
\usepackage{graphicx}

\numberwithin{equation}{section}

 3

\begin{document}

\title{\textbf{On a nonlinear heat equation associated}\\
\textbf{with Dirichlet -- Robin conditions}}
\date{}
\author{Le Thi Phuong Ngoc$^{\text{(1)}}$, Nguyen Van Y$^{\text{(3a)}}$ \and %
Alain Pham Ngoc Dinh$^{\text{(2)}}$, Nguyen Thanh Long$^{\text{(3b)}}$}
\maketitle

%\QTP{Body Math}
$^{\text{(1)}}$Nhatrang Educational College, 01 Nguyen Chanh Str., Nhatrang
City, Vietnam.\medskip

$\qquad $E-mail: ngocltp@gmail.com, ngoc1966@gmail.com\medskip

$^{\text{(2)}}$MAPMO, UMR 6628, b\^{a}t. Math\'{e}matiques, University of Orl%
\'{e}ans, BP 6759, 45067 Orl\'{e}eans Cedex 2, France.\medskip

$\qquad $E-mail: alain.pham@univ-orleans.fr, alain.pham@math.cnrs.fr\medskip

$^{\text{(3)}}$Department of Mathematics and Computer Science, University of
Natural Science, Vietnam National University HoChiMinh City, 227 Nguyen Van
Cu Str., Dist.5, HoChiMinh City, Vietnam.\medskip

$\qquad ^{\text{(3a)}}$E-mail: nguyenvanyhv@gmail.com\medskip

$\qquad ^{\text{(3b)}}$E-mail: longnt@hcmc.netnam.vn,
longnt2@gmail.com\medskip

\textbf{Abstract}. \textit{This paper is devoted to the study a nonlinear
heat equation associated with Dirichlet-Robin conditions. At first, we use
the Faedo -- Galerkin and the compactness method to prove existence and
uniqueness results}$.$\textit{\ Next, we consider the properties of
solutions. We obtain that if the initial condition is bounded then so is the
solution and we also get asymptotic behavior of solutions as }$t\rightarrow
+\infty .$\textit{\ Finally, we give numerical results}.$\bigskip $

\textbf{Keywords}: \textit{Faedo - Galerkin method, nonlinear heat equation,
Robin conditions,} \textit{Asymptotic behavior of the solution.}

\textbf{AMS subject classification}: 34B60, 35K55, 35Q72, 80A30.

\textbf{Address for correspondence}: Nguyen Thanh Long.

\section{\textbf{Introduction}}

\qquad In this paper, we consider the following nonlinear heat equation
\begin{equation}
\begin{tabular}{l}
$u_{t}-\frac{\partial }{\partial x}\left[ \mu \left( x,t\right) u_{x}\right]
+f(u)=f_{1}(x,t),\text{ }0<x<1,\text{ }0<t<T,$%
\end{tabular}
\tag{1.1}  \label{1}
\end{equation}%
associated with conditions%
\begin{equation}
\begin{tabular}{l}
$u_{x}(0,t)=h_{0}u(0,t)+g_{0}(t),\text{ }-u_{x}(1,t)=h_{1}u(1,t)+g_{1}(t),$%
\end{tabular}
\tag{1.2}  \label{2}
\end{equation}%
and initial condition%
\begin{equation}
u(x,0)=u_{0}(x),  \tag{1.3}  \label{3}
\end{equation}%
where $u_{0},$ $\mu ,$ $f,$ $f_{1},$ $g_{0},$ $g_{1}$ are given functions
satisfying conditions, which will be specified later, and $h_{0},$ $%
h_{1}\geq 0$ are given constants, with $h_{0}+h_{1}>0.$

The conditions (\ref{2})\ are commonly known as Dirichlet -- Robin
conditions. They connect Dirichlet and Neumann conditions. Theses conditions
arise from the effect of excess inert electrolytes in an electrochemical
system through perturbation analysis (\cite{2}, \cite{6}, \cite{7}, \cite{8}%
).

The governing equations (\ref{2}) are the equation usually used in a
diffusion, convection, migration transport system with electrochemical
reactions occurring at the boundary electrodes and submitted to non linear
constraints.

In electrochemistry, the oxidation-reduction reactions producing the current
is modeled by a non linear elliptic boundary value problem, linearization of
which gives the Dirichlet -- Robin conditions (\cite{3}). Theses conditions
also appear in the response of an electrochemical thin film, such as
separation in a micro -- battery. His analyze is made by solving the Poisson
-- Nernst -- Planck equation subject to boundary conditions appropriate
(Dirichlet -- Robin conditions) for an electrolytic cell (\cite{4}).

The paper consists of six sections. In Section 2, we present some
preliminaries. Using the Faedo -- Galerkin method and the compactness
method, in Section 3, we establish the existence of a unique weak solution
of the problem (\ref{1}) -- (\ref{3})\ on $(0,T),$ for every $T>0.$ In
section 4,\ we prove that if the initial condition is bounded, then so is
the solution. In section 5,\ we study asymptotic behavior of the solution as
$t\rightarrow +\infty .$ In section 6 we give numerical results.

\section{\textbf{Preliminaries}}

$\qquad $Put $\Omega =(0,1),$ $Q_{T}=\Omega \times (0,T).$ We will omit the
definitions of the usual function spaces and denote them by the notations $%
L^{p}=\,L^{p}(\Omega ),$ $H^{m}=H^{m}\left( \Omega \right) .$ Let $\langle
\cdot ,\cdot \rangle $ be either the scalar product in $L^{2}$ or the dual
pairing of a continuous linear functional and an element of a function
space. The notation $||\cdot ||$ stands for the norm in $L^{2}$ and we
denote by $||\cdot ||_{X}$ the norm in the Banach space $X.$ We call $%
X^{\prime }$ the dual space of $X.$ We denote $L^{p}(0,T;X),$ $1\leq p\leq
\infty $ the Banach space of real functions $u:(0,T)\rightarrow X$
measurable, such that $||u||_{L^{p}(0,T;X)}<+\infty ,$ with%
\begin{equation*}
||u||_{L^{p}(0,T;X)}=\left\{
\begin{tabular}{lll}
$\left( \int_{0}^{T}||u(t)||_{X}^{p}dt\right) ^{1/p},$\medskip & if & $%
1\leq p<\infty ,$\medskip \\
$\underset{0<t<T}{ess\sup }||u(t)||_{X},$ & if & $p=\infty .$%
\end{tabular}%
\right.
\end{equation*}

Let $u(t),$ $u^{\prime }(t)=u_{t}(t)=\,\overset{\cdot }{u}(t),$ $%
u_{x}(t)=\bigtriangledown u(t),$ $u_{xx}(t)=\Delta u(t),$ denote $u(x,t),$ $%
\frac{\partial u}{\partial t}(x,t),$ $\frac{\partial u}{\partial x}(x,t),$ $%
\frac{\partial ^{2}u}{\partial x^{2}}(x,t),$ respectively.

On $H^{1}$ we shall use the following norms $\left\Vert v\right\Vert
_{H^{1}}=\left( \left\Vert v\right\Vert ^{2}+\left\Vert v_{x}\right\Vert
^{2}\right) ^{1/2},$\ $\left\Vert v\right\Vert _{i}=\left(
v^{2}(i)+\left\Vert v_{x}\right\Vert ^{2}\right) ^{1/2},$ $i=0,1.$

Let $\mu \in C^{0}\left( \overline{Q}_{T}\right) ,$ with $\mu (x,t)\geq \mu
_{0}>0,$ for all $(x,t)\in \overline{Q}_{T},$ and the constants $h_{0},$ $%
h_{1}\geq 0,$ with $h_{0}+h_{1}>0,$ we consider a familly of symmetric
bilinear forms $\{a(t;\cdot ,\cdot )\}_{0\leq t\leq T}$ on $H^{1}\times
H^{1} $ as follows%
\begin{equation}
\begin{tabular}{l}
$a(t;u,v)=\int\nolimits_{0}^{1}\mu (x,t)u_{x}(x)v_{x}(x)dx+h_{0}\mu \left(
0,t\right) u(0)v(0)+h_{1}\mu \left( 1,t\right) u(1)v(1)\bigskip $ \\
$\ \ \ \ =\left\langle \mu (t)u_{x},v_{x}\right\rangle +h_{0}\mu \left(
0,t\right) u(0)v(0)+h_{1}\mu \left( 1,t\right) u(1)v(1),$ for all$\mathit{\ }%
u,v\in H^{1},$ $0\leq t\leq T.$%
\end{tabular}
\tag{2.1}  \label{b1}
\end{equation}

Then we have the following lemmas.

\textbf{Lemma 2.1}. \textit{The imbedding} $H^{1}\hookrightarrow
C^{0}([0,1]) $ \textit{is compact and}%
\begin{equation}
\left\{
\begin{tabular}{l}
$\left\Vert v\right\Vert _{C^{0}(\overline{\Omega })}\leq \sqrt{2}\left\Vert
v\right\Vert _{H^{1}},\text{ \textit{for all}}\mathit{\ }v\in H^{1},\bigskip
$ \\
$\left\Vert v\right\Vert _{C^{0}(\overline{\Omega })}\leq \sqrt{2}\left\Vert
v\right\Vert _{i},\text{ \textit{for all}}\mathit{\ }v\in H^{1},$ $i=0,1.$%
\end{tabular}%
\right.  \tag{2.2}  \label{b2}
\end{equation}

\textbf{Lemma 2.2}. \textit{Let} $\mu \in C^{0}\left( \overline{Q}%
_{T}\right) ,$ \textit{with} $\mu (x,t)\geq \mu _{0}>0,$ \textit{for all} $%
(x,t)\in \overline{Q}_{T},$ \textit{and the constants} $h_{0},$ $h_{1}\geq
0, $ with $h_{0}+h_{1}>0.$ \textit{Then, the symmetric bilinear form }$%
a(t;\cdot ,\cdot )$ \textit{is continuous on} $H^{1}\times H^{1}$ \textit{%
and coercive on} $H^{1},$ i.e.,%
\begin{equation}
\begin{tabular}{lll}
$(i)$ &  & $\left\vert a(t;u,v)\right\vert \leq a_{T}\left\Vert u\right\Vert
_{H^{1}}\left\Vert v\right\Vert _{H^{1}},\bigskip $ \\
$(ii)$ &  & $a(t;v,v)\geq a_{0}\left\Vert v\right\Vert _{H^{1}}^{2},$%
\end{tabular}%
\text{ \ \ \ \ \ \ \ \ \ \ \ \ \ \ \ \ \ \ \ \ \ \ \ \ \ \ \ \ \ \ \ \ \ }
\tag{2.3}  \label{b3}
\end{equation}%
\ \textit{for all} $u,$ $v\in H^{1},$ $0\leq t\leq T,$ \textit{where }$%
a_{T}=\left( 1+2h_{0}+2h_{1}\right) \underset{(x,t)\in \overline{Q}_{T}}{%
\sup }\mu (x,t),$ \textit{and}%
\begin{equation}
a_{0}=a_{0}(\mu _{0},h_{0},h_{1})=\left\{
\begin{tabular}{ll}
$\mu _{0}\min \{h_{0},\frac{1}{2}\},$ & $h_{0}>0,$ $h_{1}\geq 0,\bigskip $
\\
$\mu _{0}\min \{h_{1},\frac{1}{2}\},$ & $h_{1}>0,$ $h_{0}\geq 0.$%
\end{tabular}%
\right.  \tag{2.4}  \label{b4}
\end{equation}

The proofs of these lemmas are straightforward. We shall omit the details.

\textbf{Remark 2.1}. It follows from (\ref{b2}) that on $H^{1},$ $%
v\longmapsto \left\Vert v\right\Vert _{H^{1}}\ $and $v\longmapsto \left\Vert
v\right\Vert _{i}$ are two equivalent norms satisfying%
\begin{equation}
\begin{tabular}{l}
$\frac{1}{\sqrt{3}}\left\Vert v\right\Vert _{H^{1}}\leq \left\Vert
v\right\Vert _{i}\leq \sqrt{3}\left\Vert v\right\Vert _{H^{1}},\text{\ for
all}\mathit{\ }v\in H^{1},\text{ }i=0,1.$%
\end{tabular}
\tag{2.5}  \label{b5}
\end{equation}

\section{\textbf{The existence and uniqueness theorem}}

We make the following assumptions:%
\begin{equation*}
\begin{tabular}{ll}
$\left( H_{1}\right) $ & $h_{0}\geq 0$ and $h_{1}\geq 0,$ with $%
h_{0}+h_{1}>0,\bigskip $ \\
$\left( H_{2}\right) $ & $u_{0}\in L^{2},\bigskip $ \\
$\left( H_{3}\right) $ & $g_{0},$ $g_{1}\in W^{1,1}(0,T),\bigskip $ \\
$\left( H_{4}\right) $ & $\mu \in C^{1}([0,1]\times \lbrack 0,T]),$ $\mu
(x,t)\geq \mu _{0}>0,$ $\forall (x,t)\in \lbrack 0,1]\times \lbrack
0,T],\bigskip $ \\
$\left( H_{5}\right) $ & $f_{1}\in L^{1}(0,T;L^{2}),\bigskip $ \\
$\left( H_{6}\right) $ & $f\in C^{0}(%
%TCIMACRO{\U{211d} }%
%BeginExpansion
\mathbb{R}
%EndExpansion
)$ satisfies the condition, t$\text{here\thinspace \thinspace exist positive}%
\,\,\text{constants\thinspace \thinspace }C_{1},\text{ }C_{1}^{\prime },%
\text{ }C_{2}\,\text{ and }p>1,\bigskip $ \\
&
\begin{tabular}{ll}
$(i)$ & $uf(u)\geq C_{1}\left\vert u\right\vert ^{p}-C_{1}^{\prime
},\bigskip $ \\
$(2i)$ & $\left\vert f(u)\right\vert \leq C_{2}(1+\left\vert u\right\vert
^{p-1}),$ for all $u\in
%TCIMACRO{\U{211d} }%
%BeginExpansion
\mathbb{R}
%EndExpansion
.$%
\end{tabular}%
\end{tabular}%
\end{equation*}

The weak formulation of the initial boundary valued (\ref{1}) -- (\ref{3})
can then be given in the following manner: Find $u(t)$ defined in the open
set $(0,T)$ such that $u(t)$ satisfies the following variational problem%
\begin{equation}
\begin{tabular}{l}
$\frac{d}{dt}\langle u(t),v\rangle +a(t,u(t),v)+\langle f(u),v\rangle
=\langle f_{1}(t),v\rangle -\mu \left( 0,t\right) g_{0}(t)v(0)-\mu \left(
1,t\right) g_{1}(t)v(1),$%
\end{tabular}
\tag{3.1}  \label{c1}
\end{equation}%
$\forall v\in H^{1},$ and the initial condition%
\begin{equation}
\begin{tabular}{l}
$u(0)=u_{0}.$%
\end{tabular}
\tag{3.2}  \label{c2}
\end{equation}

We then have the following theorem.

\textbf{Theorem 3.1.} \textit{Let }$T>0$ \textit{and }$(H_{1})-(H_{6})$%
\textit{\ hold. Then, there exists a weak solution }$u$\textit{\ of problem }%
(\ref{1}) -- (\ref{3})\textit{\ such that}%
\begin{equation}
\left\{
\begin{tabular}{l}
$u\in L^{2}(0,T;H^{1})\cap L^{\infty }(0,T;L^{2}),\bigskip $ \\
$tu\in L^{\infty }(0,T;H^{1}),$ $tu_{t}\in L^{2}(0,T;L^{2}).$%
\end{tabular}%
\right.  \tag{3.3}  \label{c3}
\end{equation}

\textit{Furthermore, if }$f$ \textit{satisfies the following condition, in
addition,}%
\begin{equation*}
\left( H_{7}\right)
\begin{tabular}{l}
$(y-z)\left( f(y)-f(z)\right) \geq -\delta \left\vert y-z\right\vert ^{2},%
\text{ \textit{for all} }y,$ $z\in
%TCIMACRO{\U{211d} }%
%BeginExpansion
\mathbb{R}
%EndExpansion
,\text{ \textit{with} }\delta >0,$%
\end{tabular}%
\end{equation*}%
\textit{then the solution is unique}.

\textbf{Proof}. The proof consists of several steps.

\textbf{Step 1:} \textit{The Faedo -- Galerkin approximation} (introduced by
Lions \cite{5}).

Let $\{w_{j}\}$ be a denumerable base of $H^{1}.$ We find the approximate
solution of the problem (\ref{1}) -- (\ref{3}) in the form%
\begin{equation}
\begin{tabular}{l}
$u_{m}(t)=\sum_{j=1}^{m}c_{mj}(t)w_{j},$%
\end{tabular}
\tag{3.4}  \label{c4}
\end{equation}%
where the coefficients $c_{mj}$ satisfy the system of linear differential
equations%
\begin{equation}
\left\{
\begin{tabular}{l}
$\langle u_{m}^{\prime }(t),w_{j}\rangle +a(t;u_{m}(t),w_{j})+\langle
f(u_{m}(t)),w_{j}\rangle \bigskip $ \\
$\ \ \ \ \ \ \ \ \ \ \ \ \ \ =\langle f_{1}(t),w_{j}\rangle -\mu \left(
0,t\right) g_{0}(t)w_{j}(0)-\mu \left( 1,t\right) g_{1}(t)w_{j}(1),$ $1\leq
j\leq m,\bigskip $ \\
$u_{m}(0)=u_{0m},$%
\end{tabular}%
\right.  \tag{3.5}  \label{c5}
\end{equation}%
where%
\begin{equation}
\begin{tabular}{l}
$u_{0m}=\sum_{j=1}^{m}\alpha _{mj}w_{j}\rightarrow u_{0}$ strongly in $%
L^{2}. $%
\end{tabular}
\tag{3.6}  \label{c6}
\end{equation}

It is clear that for each $m$ there exists a solution $u_{m}(t)$ in form (%
\ref{c4}) which satisfies (\ref{c5}) and (\ref{c6}) almost everywhere on $%
0\leq t\leq T_{m}$ for some $T_{m},$ $0<T_{m}\leq T.$ The following
estimates allow one to take $T_{m}=T$\ for all $m.$

\textit{Step 2. A priori estimates}.

a) \textit{The first estimate}. Multiplying the $j^{th}$\ equation of (\ref%
{c5})\ by $c_{mj}(t)$\ and summing up with respect to $j,$ afterwards,
integrating by parts with respect to the time variable from $0$ to $t,$ we
get after some rearrangements%
\begin{equation}
\begin{tabular}{l}
$\left\Vert u_{m}(t)\right\Vert
^{2}+2\int\nolimits_{0}^{t}a(s;u_{m}(s),u_{m}(s))ds+2\int\nolimits_{0}^{t}%
\langle f(u_{m}(s)),u_{m}(s)\rangle ds\bigskip $ \\
$\ \ \ \ \ \ \ \ \ \ \ \ \ \ \ \ \ =\left\Vert u_{0m}\right\Vert
^{2}+2\int\nolimits_{0}^{t}\langle f_{1}(s),u_{m}(s)\rangle ds\bigskip $ \\
$\ \ \ \ \ \ \ \ \ \ \ \ \ \ \ \ \ \ -2\int\nolimits_{0}^{t}\mu \left(
0,s\right) g_{0}(s)u_{m}(0,s)ds-2\int\nolimits_{0}^{t}\mu \left( 1,s\right)
g_{1}(s)u_{m}(1,s)ds.$%
\end{tabular}
\tag{3.7}  \label{c7}
\end{equation}

By $u_{0m}\rightarrow u_{0}$ strongly in $L^{2},$ we have%
\begin{equation}
\begin{tabular}{l}
$\left\Vert u_{0m}\right\Vert ^{2}\leq C_{0},$ \ for all \ $m,$%
\end{tabular}
\tag{3.8}  \label{c8}
\end{equation}%
where $C_{0}$ always indicates a bound depending on $u_{0}.$

By the assumptions $(H_{6},(i)),$ and using the inequalities (\ref{b2}), (%
\ref{b3}), and with $\beta >0,$ we estimate without difficulty the following
terms in (\ref{c7}) as follows%
\begin{equation}
\begin{tabular}{l}
$2\int\nolimits_{0}^{t}a(s;u_{m}(s),u_{m}(s))ds\geq
2a_{0}\int\nolimits_{0}^{t}\left\Vert u_{m}(s)\right\Vert _{H^{1}}^{2}ds,$%
\end{tabular}
\tag{3.9}  \label{c9}
\end{equation}%
\begin{equation}
\begin{tabular}{l}
$2\int\nolimits_{0}^{t}\langle f(u_{m}(s)),u_{m}(s)\rangle ds\geq
2C_{1}\int\nolimits_{0}^{t}\left\Vert u_{m}(s)\right\Vert
_{L^{p}}^{p}ds-2TC_{1}^{\prime },$%
\end{tabular}
\tag{3.10}  \label{c10}
\end{equation}%
\begin{equation}
\begin{tabular}{l}
$2\int\nolimits_{0}^{t}\langle f_{1}(s),u_{m}(s)\rangle ds\leq \left\Vert
f_{1}\right\Vert _{L^{1}(0,T;L^{2})}+\int\nolimits_{0}^{t}\left\Vert
f_{1}(s)\right\Vert \left\Vert u_{m}(s)\right\Vert ^{2}ds,$%
\end{tabular}
\tag{3.11}  \label{c11}
\end{equation}%
\begin{equation}
\begin{tabular}{l}
$-2\int\nolimits_{0}^{t}\mu \left( 0,s\right) g_{0}(s)u_{m}(0,s)ds\leq 2%
\sqrt{2}\left\Vert \mu \right\Vert _{L^{\infty }(Q_{T})}\left\Vert
g_{0}\right\Vert _{L^{\infty }}\int\nolimits_{0}^{t}\left\Vert
u_{m}(s)\right\Vert _{H^{1}}ds\bigskip $ \\
$\ \ \ \ \ \ \ \ \ \ \ \ \ \ \ \ \ \ \ \ \ \ \ \ \ \ \ \ \ \ \ \ \ \ \ \ \ \
\ \ \ \leq \frac{2}{\beta }T\left\Vert \mu \right\Vert _{L^{\infty
}(Q_{T})}^{2}\left\Vert g_{0}\right\Vert _{L^{\infty }}^{2}+\beta
\int\nolimits_{0}^{t}\left\Vert u_{m}(s)\right\Vert _{H^{1}}^{2}ds,$%
\end{tabular}
\tag{3.12}  \label{c12}
\end{equation}%
\begin{equation}
\begin{tabular}{l}
$-2\int\nolimits_{0}^{t}\mu \left( 1,s\right) g_{1}(s)u_{m}(1,s)ds\leq 2%
\sqrt{2}\left\Vert \mu \right\Vert _{L^{\infty }(Q_{T})}\left\Vert
g_{1}\right\Vert _{L^{\infty }}\int\nolimits_{0}^{t}\left\Vert
u_{m}(s)\right\Vert _{H^{1}}ds\bigskip $ \\
$\ \ \ \ \ \ \ \ \ \ \ \ \ \ \ \ \ \ \ \ \ \ \ \ \ \ \ \ \ \ \ \ \ \ \ \ \ \
\ \ \ \leq \frac{2}{\beta }T\left\Vert \mu \right\Vert _{L^{\infty
}(Q_{T})}^{2}\left\Vert g_{1}\right\Vert _{L^{\infty }}^{2}+\beta
\int\nolimits_{0}^{t}\left\Vert u_{m}(s)\right\Vert _{H^{1}}^{2}ds,$%
\end{tabular}
\tag{3.13}  \label{c13}
\end{equation}%
for all $\beta >0.$ Hence, it follows from (\ref{c7}) -- (\ref{c13}) that%
\begin{equation}
\begin{tabular}{l}
$\left\Vert u_{m}(t)\right\Vert ^{2}+2(a_{0}-\beta
)\int\nolimits_{0}^{t}\left\Vert u_{m}(s)\right\Vert
_{H^{1}}^{2}ds+2C_{1}\int\nolimits_{0}^{t}\left\Vert u_{m}(s)\right\Vert
_{L^{p}}^{p}ds\bigskip $ \\
$\ \ \ \ \ \ \ \ \ \ \ \ \ \ \ \ \ \leq C_{0}+2TC_{1}^{\prime }+\left\Vert
f_{1}\right\Vert _{L^{1}(0,T;L^{2})}+\int\nolimits_{0}^{t}\left\Vert
f_{1}(s)\right\Vert \left\Vert u_{m}(s)\right\Vert ^{2}ds\bigskip $ \\
$\ \ \ \ \ \ \ \ \ \ \ \ \ \ \ \ \ \ \ +\frac{2}{\beta }T\left\Vert \mu
\right\Vert _{L^{\infty }(Q_{T})}^{2}\left( \left\Vert g_{0}\right\Vert
_{L^{\infty }}^{2}+\left\Vert g_{1}\right\Vert _{L^{\infty }}^{2}\right) .$%
\end{tabular}
\tag{3.14}  \label{c14}
\end{equation}

Choosing $\beta =\frac{1}{2}a_{0},$ we deduce from (\ref{c14}), that%
\begin{equation}
\begin{tabular}{l}
$S_{m}(t)\leq C_{T}^{(1)}+\int\nolimits_{0}^{t}C_{T}^{(2)}(s)S_{m}(s)ds,$%
\end{tabular}
\tag{3.15}  \label{c15}
\end{equation}%
where%
\begin{equation}
\left\{
\begin{tabular}{l}
$S_{m}(t)=\left\Vert u_{m}(t)\right\Vert
^{2}+a_{0}\int\nolimits_{0}^{t}\left\Vert u_{m}(s)\right\Vert
_{H^{1}}^{2}ds+2C_{1}\int\nolimits_{0}^{t}\left\Vert u_{m}(s)\right\Vert
_{L^{p}}^{p}ds,\bigskip $ \\
$C_{T}^{(1)}=C_{0}+2TC_{1}^{\prime }+\left\Vert f_{1}\right\Vert
_{L^{1}(0,T;L^{2})}+\frac{4}{a_{0}}T\left\Vert \mu \right\Vert _{L^{\infty
}(Q_{T})}^{2}\left( \left\Vert g_{0}\right\Vert _{L^{\infty
}}^{2}+\left\Vert g_{1}\right\Vert _{L^{\infty }}^{2}\right) ,\bigskip $ \\
$C_{T}^{(2)}(s)=\left\Vert f_{1}(s)\right\Vert ,$ $\ C_{T}^{(2)}\in
L^{1}(0,T).$%
\end{tabular}%
\right.  \tag{3.16}  \label{c16}
\end{equation}

By the Gronwall's lemma, we obtain from (\ref{c15}), that%
\begin{equation}
\begin{tabular}{l}
$S_{m}(t)\leq C_{T}^{(1)}\exp \left(
\int\nolimits_{0}^{t}C_{T}^{(2)}(s)ds\right) \leq C_{T},$%
\end{tabular}
\tag{3.17}  \label{c17}
\end{equation}%
for all $m\in
%TCIMACRO{\U{2115} }%
%BeginExpansion
\mathbb{N}
%EndExpansion
,$ for all $t,$\ $0\leq t\leq T_{m}\leq T,$ \ i.e., $T_{m}=T,$ where $C_{T}$
always indicates a bound depending on $T.$

b) \textit{The second estimate}. Multiplying the $j^{th}$ equation of the
system (\ref{c5}) by $t^{2}c_{mj}^{\prime }(t)$ and summing up with respect
to $j,$ we have%
\begin{equation}
\begin{tabular}{l}
$\left\Vert tu_{m}^{\prime }(t)\right\Vert
^{2}+t^{2}a(t;u_{m}(t),u_{m}^{\prime }(t))+\langle
tf(u_{m}(t)),tu_{m}^{\prime }(t)\rangle \bigskip $ \\
$\ \ \ \ \ \ \ \ \ \ \ \ \ =\langle tf_{1}(t),tu_{m}^{\prime }(t)\rangle
-t^{2}\mu \left( 0,t\right) g_{0}(t)u_{m}^{\prime }(0,t)-t^{2}\mu \left(
1,t\right) g_{1}(t)u_{m}^{\prime }(1,t).$%
\end{tabular}
\tag{3.18}  \label{c18}
\end{equation}

First, we need the following lemmas.

\textbf{Lemma 3.2}.%
\begin{equation}
\begin{tabular}{l}
(i) $\ \ \frac{\partial a}{\partial t}(t;u,v)=\langle \mu ^{\prime }\left(
\cdot ,t\right) u_{x},v_{x}\rangle +h_{0}\mu ^{\prime }\left( 0,t\right)
u(0)v(0)+h_{1}\mu ^{\prime }\left( 1,t\right) u(1)v(1),$ \textit{for all} $%
u,v\in H^{1},\bigskip $ \\
(ii) $\ \ \left\vert \frac{\partial a}{\partial t}(t;u,v)\right\vert \leq
\widetilde{a}_{T}\left\Vert u\right\Vert _{H^{1}}\left\Vert v\right\Vert
_{H^{1}},\text{\ \textit{for all} }u,v\in H^{1},\bigskip $ \\
(iii) $\ \ \frac{d}{dt}a(t;u_{m}(t),u_{m}(t))=2a(t;u_{m}(t),u_{m}^{\prime
}(t))+\frac{\partial a}{\partial t}(t;u_{m}(t),u_{m}(t)),$%
\end{tabular}
\tag{3.19}  \label{c19}
\end{equation}%
\textit{where} $\widetilde{a}_{T}=\left( 1+2h_{0}+2h_{1}\right) \underset{%
(x,t)\in \lbrack 0,1]\times \lbrack 0,T]}{\sup }\mu ^{\prime }\left(
x,t\right) .\bigskip $

\textbf{Lemma 3.3}. \textit{Put} $\lambda _{0}=\left( \frac{C_{1}^{\prime }}{%
C_{1}}\right) ^{1/p},$ $m_{0}=\int\nolimits_{-\lambda _{0}}^{\lambda
_{0}}\left\vert f(y)\right\vert dy,$\ and \ $\overline{f}(z)=\int%
\nolimits_{0}^{z}f(y)dy,$ $z\in
%TCIMACRO{\U{211d} }%
%BeginExpansion
\mathbb{R}
%EndExpansion
.$

$\qquad $\textit{Then we have}%
\begin{equation}
\begin{tabular}{l}
$-m_{0}\leq \overline{f}(z)\leq C_{2}(\left\vert z\right\vert +\frac{1}{p}%
\left\vert z\right\vert ^{p}),$\ $\ \forall z\in
%TCIMACRO{\U{211d} }%
%BeginExpansion
\mathbb{R}
%EndExpansion
.$%
\end{tabular}
\tag{3.20}  \label{c20}
\end{equation}

The proofs of these lemmas are straightforward. We shall omit the details.$%
\blacksquare $

By (\ref{c19})$_{3}$, we rewrite (\ref{c18}) as follows%
\begin{equation}
\begin{tabular}{l}
$2\left\Vert tu_{m}^{\prime }(t)\right\Vert ^{2}+\frac{d}{dt}%
a(t;tu_{m}(t),tu_{m}(t))+2\langle tf(u_{m}(t)),tu_{m}^{\prime }(t)\rangle
\bigskip $ \\
$\ \ \ \ \ \ \ \ \ \ \ \ \ \ =2ta(t;u_{m}(t),u_{m}(t))+\frac{\partial a}{%
\partial t}(t;tu_{m}(t),tu_{m}(t))+2\langle tf_{1}(t),tu_{m}^{\prime
}(t)\rangle \bigskip $ \\
$\ \ \ \ \ \ \ \ \ \ \ \ \ \ \ -2t^{2}\mu \left( 0,t\right)
g_{0}(t)u_{m}^{\prime }(0,t)-2t^{2}\mu \left( 1,t\right)
g_{1}(t)u_{m}^{\prime }(1,t).$%
\end{tabular}
\tag{3.21}  \label{c21}
\end{equation}

Integrating (\ref{c21}), we get%
\begin{equation}
\begin{tabular}{l}
$2\int\nolimits_{0}^{t}\left\Vert su_{m}^{\prime }(s)\right\Vert
^{2}ds+a(t;tu_{m}(t),tu_{m}(t))+2\int\nolimits_{0}^{t}\langle
sf(u_{m}(s)),su_{m}^{\prime }(s)\rangle ds\bigskip $ \\
$\ \ \ \ \ \ \ \ \ \ \ \ \
=2\int\nolimits_{0}^{t}sa(s;u_{m}(s),u_{m}(s))ds+\int\nolimits_{0}^{t}%
\frac{\partial a}{\partial t}(s;su_{m}(s),su_{m}(s))ds+2\int%
\nolimits_{0}^{t}\langle sf_{1}(s),su_{m}^{\prime }(s)\rangle ds\bigskip $
\\
$\ \ \ \ \ \ \ \ \ \ \ \ \ -2\int\nolimits_{0}^{t}s^{2}\mu \left(
0,s\right) g_{0}(s)u_{m}^{\prime }(0,s)ds-2\int\nolimits_{0}^{t}s^{2}\mu
\left( 1,s\right) g_{1}(s)u_{m}^{\prime }(1,s)ds.$%
\end{tabular}
\tag{3.22}  \label{c22}
\end{equation}

We shall estimate the terms of (\ref{c22}) as follows.%
\begin{equation}
\begin{tabular}{l}
$a(t;tu_{m}(t),tu_{m}(t))\geq a_{0}\left\Vert tu_{m}(t)\right\Vert
_{H^{1}}^{2},$%
\end{tabular}
\tag{3.23}  \label{c23}
\end{equation}%
\begin{equation}
\begin{tabular}{l}
$2\int\nolimits_{0}^{t}\langle sf(u_{m}(s)),su_{m}^{\prime }(s)\rangle
ds=2\int\nolimits_{0}^{t}s^{2}ds\frac{d}{ds}\int\nolimits_{0}^{1}dx\int%
\nolimits_{0}^{u_{m}(x,s)}f(y)dy\bigskip $ \\
$\ \ \ \ \ \ \ \ \ \ \ \ \ \ \ \ \ \ \ \ \ \ \ \ \ \ \ \ \ \ \ \ \ \ \
=2\int\nolimits_{0}^{t}s^{2}ds\frac{d}{ds}\int\nolimits_{0}^{1}\overline{f}%
(u_{m}(x,s))dx\bigskip $ \\
$\ \ \ \ \ \ \ \ \ \ \ \ \ \ \ \ \ \ \ \ \ \ \ \ \ \ \ \ \ \ \ \ \ \ \
=2\int\nolimits_{0}^{t}\left[ \frac{d}{ds}\left( s^{2}\int\nolimits_{0}^{1}%
\overline{f}(u_{m}(x,s))dx\right) -2s\int\nolimits_{0}^{1}\overline{f}%
(u_{m}(x,s))dx\right] ds\bigskip $ \\
$\ \ \ \ \ \ \ \ \ \ \ \ \ \ \ \ \ \ \ \ \ \ \ \ \ \ \ \ \ \ \ \ \ \ \
=2t^{2}\int\nolimits_{0}^{1}\overline{f}(u_{m}(x,t))dx-4\int%
\nolimits_{0}^{t}sds\int\nolimits_{0}^{1}\overline{f}(u_{m}(x,s))dx\bigskip
$ \\
$\ \ \ \ \ \ \ \ \ \ \ \ \ \ \ \ \ \ \ \ \ \ \ \ \ \ \ \ \ \ \ \ \ \ \ \geq
-2T^{2}m_{0}-4C_{2}\int\nolimits_{0}^{t}s\left[ \left\Vert
u_{m}(s)\right\Vert _{L^{1}}+\frac{1}{p}\left\Vert u_{m}(s)\right\Vert
_{L^{p}}^{p}\text{ }\right] ds\bigskip $ \\
$\ \ \ \ \ \ \ \ \ \ \ \ \ \ \ \ \ \ \ \ \ \ \ \ \ \ \ \ \ \ \ \ \ \ \ \geq
-2T^{2}m_{0}-4TC_{2}\left[ T\left\Vert u_{m}\right\Vert _{L^{\infty
}(0,T;L^{2})}+\frac{1}{p}\frac{1}{2C_{1}}S_{m}(t)\right] \geq -C_{T},$%
\end{tabular}
\tag{3.24}  \label{c24}
\end{equation}

\begin{equation}
\begin{tabular}{l}
$2\int\nolimits_{0}^{t}sa(s;u_{m}(s),u_{m}(s))ds\leq
2Ta_{T}\int\nolimits_{0}^{t}\left\Vert u_{m}(s)\right\Vert
_{H^{1}}^{2}ds\leq 2Ta_{T}\frac{1}{a_{0}}S_{m}(t)\leq C_{T},$%
\end{tabular}
\tag{3.25}  \label{c25}
\end{equation}%
\begin{equation}
\begin{tabular}{l}
$\int\nolimits_{0}^{t}\frac{\partial a}{\partial t}%
(s;su_{m}(s),su_{m}(s))ds\leq \widetilde{a}_{T}\int\nolimits_{0}^{t}\left%
\Vert su_{m}(s)\right\Vert _{H^{1}}^{2}ds\leq T^{2}\widetilde{a}%
_{T}\int\nolimits_{0}^{t}\left\Vert u_{m}(s)\right\Vert
_{H^{1}}^{2}ds\bigskip $ \\
\ \ \ \ \ \ \ \ \ \ \ \ \ \ \ \ \ \ \ \ \ \ \ \ \ \ \ \ \ \ \ \ \ \ \ \ \ \
\ \ $\leq T^{2}\widetilde{a}_{T}\frac{1}{a_{0}}S_{m}(t)\leq C_{T},$%
\end{tabular}
\tag{3.26}  \label{c26}
\end{equation}%
\begin{equation}
\begin{tabular}{l}
$2\int\nolimits_{0}^{t}\langle sf_{1}(s),su_{m}^{\prime }(s)\rangle ds\leq
2\int\nolimits_{0}^{t}\left\Vert sf_{1}(s)\right\Vert \left\Vert
su_{m}^{\prime }(s)\right\Vert ds\leq \int\nolimits_{0}^{t}\left\Vert
sf_{1}(s)\right\Vert ^{2}ds+\int\nolimits_{0}^{t}\left\Vert su_{m}^{\prime
}(s)\right\Vert ^{2}ds\bigskip $ \\
\ \ \ \ \ \ \ \ \ \ \ \ \ \ \ \ \ \ \ \ \ \ \ \ \ \ \ \ \ \ \ \ $\leq
T^{2}\int\nolimits_{0}^{T}\left\Vert f_{1}(s)\right\Vert
^{2}ds+\int\nolimits_{0}^{t}\left\Vert su_{m}^{\prime }(s)\right\Vert
^{2}ds\bigskip $ \\
$\ \ \ \ \ \ \ \ \ \ \ \ \ \ \ \ \ \ \ \ \ \ \ \ \ \ \ \ \ \ \ \leq
C_{T}+\int\nolimits_{0}^{t}\left\Vert su_{m}^{\prime }(s)\right\Vert
^{2}ds. $%
\end{tabular}
\tag{3.27}  \label{c27}
\end{equation}

By using integration by parts, it follows that%
\begin{equation}
\begin{tabular}{l}
$\left\vert -2\int\nolimits_{0}^{t}s^{2}\mu \left( 0,s\right)
g_{0}(s)u_{m}^{\prime }(0,s)ds\right\vert \bigskip $ \\
$\ \ \ \ \ \ \ \ =\left\vert -2t^{2}\mu \left( 0,t\right)
g_{0}(t)u_{m}(0,t)+2\int\nolimits_{0}^{t}\left[ s^{2}\mu \left( 0,s\right)
g_{0}(s)\right] ^{\prime }u_{m}(0,s)ds\right\vert \bigskip $ \\
\ \ \ \ \ \ \ \ \ $\leq 2\sqrt{2}t^{2}\left\Vert \mu \right\Vert _{L^{\infty
}(Q_{T})}\left\Vert g_{0}\right\Vert _{L^{\infty }}\left\Vert
u_{m}(t)\right\Vert _{H^{1}}+2\sqrt{2}\int\nolimits_{0}^{t}\left\vert \left[
s^{2}\mu \left( 0,s\right) g_{0}(s)\right] ^{\prime }\right\vert \left\Vert
u_{m}(s)\right\Vert _{H^{1}}ds\bigskip $ \\
\ \ \ \ \ \ \ \ \ $\leq \frac{2}{\beta }T^{2}\left\Vert \mu \right\Vert
_{L^{\infty }(Q_{T})}^{2}\left\Vert g_{0}\right\Vert _{L^{\infty
}}^{2}+\beta \left\Vert tu_{m}(t)\right\Vert _{H^{1}}^{2}+2\sqrt{2}%
\int\nolimits_{0}^{t}\left\vert \left[ s^{2}\mu \left( 0,s\right) g_{0}(s)%
\right] ^{\prime }\right\vert \left\Vert u_{m}(s)\right\Vert
_{H^{1}}ds\bigskip $ \\
\ \ \ \ \ \ \ \ \ $\leq \frac{1}{\beta }C_{T}+\beta \left\Vert
tu_{m}(t)\right\Vert _{H^{1}}^{2}+2\sqrt{2}\int\nolimits_{0}^{t}\left\vert %
\left[ s^{2}\mu \left( 0,s\right) g_{0}(s)\right] ^{\prime }\right\vert
\left\Vert u_{m}(s)\right\Vert _{H^{1}}ds.$%
\end{tabular}
\tag{3.28}  \label{c28}
\end{equation}

On the other hand%
\begin{equation}
\begin{tabular}{l}
$\left\vert \left[ s^{2}\mu \left( 0,s\right) g_{0}(s)\right] ^{\prime
}\right\vert =\left\vert 2s\mu \left( 0,s\right) g_{0}(s)+s^{2}\left[ \mu
^{\prime }\left( 0,s\right) g_{0}(s)+\mu \left( 0,s\right) g_{0}^{\prime }(s)%
\right] \right\vert \bigskip $ \\
\ \ \ \ \ \ \ \ \ \ \ \ \ \ \ \ \ \ \ \ \ \ \ $\leq 2s\left\Vert \mu
\right\Vert _{L^{\infty }(Q_{T})}\left\Vert g_{0}\right\Vert _{L^{\infty
}}+s^{2}\left\Vert \mu \right\Vert _{C^{1}(\overline{Q}_{T})}\left[
\left\Vert g_{0}\right\Vert _{L^{\infty }}+\left\vert g_{0}^{\prime
}(s)\right\vert \right] \bigskip $ \\
\ \ \ \ \ \ \ \ \ \ \ \ \ \ \ \ \ \ \ \ \ \ \ $\leq s\left\Vert \mu
\right\Vert _{C^{1}(\overline{Q}_{T})}\left[ (2+T)\left\Vert
g_{0}\right\Vert _{L^{\infty }}+T\left\vert g_{0}^{\prime }(s)\right\vert %
\right] \leq sC_{T}\psi _{0}(s),$%
\end{tabular}
\tag{3.29}  \label{c29}
\end{equation}%
where%
\begin{equation}
\begin{tabular}{l}
$C_{T}=\left\Vert \mu \right\Vert _{C^{1}(\overline{Q}_{T})}\left[
(2+T)\left\Vert g_{0}\right\Vert _{L^{\infty }}+T\text{ }\right] ,$ $\ \psi
_{0}(s)=1+\left\vert g_{0}^{\prime }(s)\right\vert ,$ $\ \psi _{0}\in
L^{1}(0,T).$%
\end{tabular}
\tag{3.30}  \label{c30}
\end{equation}

Hence, we deduce from (\ref{c28}), (\ref{c29}), that%
\begin{equation}
\begin{tabular}{l}
$\left\vert -2\int\nolimits_{0}^{t}s^{2}\mu \left( 0,s\right)
g_{0}(s)u_{m}^{\prime }(0,s)ds\right\vert \leq \frac{1}{\beta }C_{T}+\beta
\left\Vert tu_{m}(t)\right\Vert _{H^{1}}^{2}+2\sqrt{2}C_{T}\int%
\nolimits_{0}^{t}\psi _{0}(s)\left\Vert su_{m}(s)\right\Vert
_{H^{1}}ds\bigskip $ \\
\ \ \ \ \ \ \ \ \ \ \ \ \ \ \ \ \ \ \ \ \ \ \ \ \ \ \ \ \ \ \ \ $\leq \frac{1%
}{\beta }C_{T}+\beta \left\Vert tu_{m}(t)\right\Vert
_{H^{1}}^{2}+2C_{T}^{2}\int\nolimits_{0}^{T}\psi
_{0}(s)ds+\int\nolimits_{0}^{t}\psi _{0}(s)\left\Vert su_{m}(s)\right\Vert
_{H^{1}}^{2}ds\bigskip $ \\
\ \ \ \ \ \ \ \ \ \ \ \ \ \ \ \ \ \ \ \ \ \ \ \ \ \ \ \ \ \ \ \ \ $\leq (1+%
\frac{1}{\beta })C_{T}+\beta \left\Vert tu_{m}(t)\right\Vert
_{H^{1}}^{2}+\int\nolimits_{0}^{t}\psi _{0}(s)\left\Vert
su_{m}(s)\right\Vert _{H^{1}}^{2}ds,$%
\end{tabular}
\tag{3.31}  \label{c31}
\end{equation}%
for all $\beta >0.$

Similarly%
\begin{equation}
\begin{tabular}{l}
$-2\int\nolimits_{0}^{t}s^{2}\mu \left( 1,s\right) g_{1}(s)u_{m}^{\prime
}(1,s)ds\leq (1+\frac{1}{\beta })C_{T}+\beta \left\Vert tu_{m}(t)\right\Vert
_{H^{1}}^{2}+\int\nolimits_{0}^{t}\psi _{1}(s)\left\Vert
su_{m}(s)\right\Vert _{H^{1}}^{2}ds,$%
\end{tabular}
\tag{3.32}  \label{c32}
\end{equation}%
for all $\beta >0,$ where%
\begin{equation}
\begin{tabular}{l}
$C_{T}=\left\Vert \mu \right\Vert _{C^{1}(\overline{Q}_{T})}\left[
(2+T)\left\Vert g_{1}\right\Vert _{L^{\infty }}+T\text{ }\right] ,$ $\ \psi
_{1}(s)=1+\left\vert g_{1}^{\prime }(s)\right\vert ,$ $\ \psi _{1}\in
L^{1}(0,T).$%
\end{tabular}
\tag{3.33}  \label{c33}
\end{equation}

It follows from (\ref{c22}) -- (\ref{c27}), (\ref{c31})\ and (\ref{c32}),
that%
\begin{equation}
\begin{tabular}{l}
$\int\nolimits_{0}^{t}\left\Vert su_{m}^{\prime }(s)\right\Vert
^{2}ds+a_{0}\left\Vert tu_{m}(t)\right\Vert _{H^{1}}^{2}\bigskip $ \\
$\ \ \ \ \ \ \ \ \ \ \ \ \ \leq (6+\frac{2}{\beta })C_{T}+2\beta \left\Vert
tu_{m}(t)\right\Vert _{H^{1}}^{2}+\int\nolimits_{0}^{t}\psi
_{0}(s)\left\Vert su_{m}(s)\right\Vert _{H^{1}}^{2}ds\bigskip $ \\
$\ \ \ \ \ \ \ \ \ \ \ \ \ \ +\int\nolimits_{0}^{t}\psi _{1}(s)\left\Vert
su_{m}(s)\right\Vert _{H^{1}}^{2}ds.$%
\end{tabular}
\tag{3.34}  \label{c34}
\end{equation}

Choosing $2\beta =\frac{1}{2}a_{0},$ we deduce from (\ref{c34}), that%
\begin{equation}
\begin{tabular}{l}
$X_{m}(t)\leq \overline{C}_{T}^{(1)}+\int\nolimits_{0}^{t}\overline{C}%
_{T}^{(2)}(s)X_{m}(s)ds,$%
\end{tabular}
\tag{3.35}  \label{c35}
\end{equation}%
where%
\begin{equation}
\left\{
\begin{tabular}{l}
$X_{m}(t)=\left\Vert tu_{m}(t)\right\Vert
_{H^{1}}^{2}+\int\nolimits_{0}^{t}\left\Vert su_{m}^{\prime }(s)\right\Vert
^{2}ds,\bigskip $ \\
$\overline{C}_{T}^{(1)}=\left( 1+\frac{2}{a_{0}}\right) (6+\frac{8}{a_{0}}%
)C_{T},\bigskip $ \\
$\overline{C}_{T}^{(2)}(s)=\left( 1+\frac{2}{a_{0}}\right) \left( \psi
_{0}(s)+\psi _{1}(s)\right) ,\ \ \ \overline{C}_{T}^{(2)}\in L^{1}(0,T).$%
\end{tabular}%
\right.  \tag{3.36}  \label{c36}
\end{equation}

By the Gronwall's lemma, we obtain from (\ref{c35}), that%
\begin{equation}
\begin{tabular}{l}
$\left\Vert tu_{m}(t)\right\Vert
_{H^{1}}^{2}+\int\nolimits_{0}^{t}\left\Vert su_{m}^{\prime }(s)\right\Vert
^{2}ds\leq \overline{C}_{T}^{(1)}\exp \left( \int\nolimits_{0}^{T}\overline{%
C}_{T}^{(2)}(s)ds\right) \leq C_{T},$%
\end{tabular}
\tag{3.37}  \label{c37}
\end{equation}%
for all $m\in
%TCIMACRO{\U{2115} }%
%BeginExpansion
\mathbb{N}
%EndExpansion
,$ for all $t\in \lbrack 0,T],$ $\forall T>0,$ where $C_{T}$ always
indicates a bound depending on $T.$

\textit{Step 3}. \textit{The limiting process}.

By (\ref{c16}), (\ref{c17}) and (\ref{c37}) we deduce that, there exists a
subsequence of $\{u_{m}\},$ still denoted by $\{u_{m}\}$ such that%
\begin{equation}
\left\{
\begin{tabular}{lll}
$u_{m}\rightarrow u$ & $\text{in}$ & $L^{\infty }(0,T;L^{2})\text{ \ weak*,}$%
\medskip \\
$u_{m}\rightarrow u$ & $\text{in}$ & $L^{2}(0,T;H^{1})\text{ \ weak,}$%
\medskip \\
$tu_{m}\rightarrow tu$ & $\text{in}$ & $L^{\infty }(0,T;H^{1})\text{ \ weak*,%
}$\medskip \\
$(tu_{m})^{\prime }\rightarrow (tu)^{\prime }$ & $\text{in}$ & $L^{2}(Q_{T})%
\text{ \ weak,}$\medskip \\
$u_{m}\rightarrow u$ & $\text{in}$ & $L^{p}(Q_{T})\text{ \ weak.}$%
\end{tabular}%
\right.  \tag{3.38}  \label{c38}
\end{equation}

Using a compactness lemma (\cite{5}, Lions, p. 57) applied to (\ref{c38})$%
_{3,4}$, we can extract from the sequence $\{u_{m}\}$ a subsequence still
denotes by $\{u_{m}\},$ such that%
\begin{equation}
\begin{tabular}{l}
$tu_{m}\rightarrow tu\text{\ \ strongly in }L^{2}(Q_{T}).$%
\end{tabular}
\tag{3.39}  \label{c39}
\end{equation}

By the Riesz- Fischer theorem, we can extract from $\{u_{m}\}$ a subsequence
still denoted by $\{u_{m}\},$ such that%
\begin{equation}
\begin{tabular}{l}
$u_{m}(x,t)\rightarrow u(x,t)\text{ \ a.e. \ }(x,t)\ \ \text{in \ }%
Q_{T}=(0,1)\times (0,T).$%
\end{tabular}
\tag{3.40}  \label{c40}
\end{equation}

Because $f$ is continuous, then%
\begin{equation}
\begin{tabular}{l}
$f(u_{m}(x,t))\rightarrow f(u(x,t))\text{ \ a.e. \ }(x,t)\ \ \text{in \ }%
Q_{T}=(0,1)\times (0,T).$%
\end{tabular}
\tag{3.41}  \label{c41}
\end{equation}

On the other hand, by $(H_{6},$ $ii),$ it follows from (\ref{c16}), (\ref%
{c17}) that%
\begin{equation}
\begin{tabular}{l}
$\left\Vert f(u_{m})\right\Vert _{L^{p^{\prime }}(Q_{T})}\leq C_{T},$%
\end{tabular}
\tag{3.42}  \label{c42}
\end{equation}%
where $C_{T}$ is a constant independent of\ $m.$

We shall now require the following lemma, the proof of which can be found in
\cite{5}.

\textbf{Lemma 3.4}. \textit{Let} $Q$ \textit{be a bounded open set of} $\
%TCIMACRO{\U{211d} }%
%BeginExpansion
\mathbb{R}
%EndExpansion
^{N}$ \textit{and} $G_{m},$ $G\in L^{q}(Q),$ $1<q<\infty ,$ \textit{such that%
},%
\begin{equation*}
\begin{tabular}{l}
$\left\Vert G_{m}\right\Vert _{L^{q}(Q)}\leq C,\text{ \textit{where} }C\text{
\textit{is a constant independent of}\ }m,$%
\end{tabular}%
\end{equation*}%
\textit{and}%
\begin{equation*}
\begin{tabular}{l}
$G_{m}\rightarrow G\text{ \ a.e. \ }(x,t)\ \ \text{in \ }Q.$%
\end{tabular}%
\end{equation*}

\textit{Then}%
\begin{equation*}
\begin{tabular}{l}
$G_{m}\rightarrow G$ $\ $\textit{in } $L^{q}(Q)$ $\ $\textit{weakly.}$%
\blacksquare $%
\end{tabular}%
\end{equation*}

Applying Lemma 3.4 with $N=2,$ $q=p^{\prime },\ G_{m}=f(u_{m}),$ $G=f(u),$
we deduce from (\ref{c41}), (\ref{c42}), that%
\begin{equation}
\begin{tabular}{l}
$f(u_{m})\rightarrow f(u)$ $\ $\textit{in } $L^{p^{\prime }}(Q_{T})$ $\ $%
\textit{weakly.}%
\end{tabular}
\tag{3.43}  \label{c43}
\end{equation}

Passing to the limit in (\ref{c5}) by (\ref{c6}), (\ref{c38}), (\ref{c43}),
we have satisfying the equation%
\begin{equation}
\left\{
\begin{tabular}{l}
$\frac{d}{dt}\langle u(t),v\rangle +a(t,u(t),v)+\langle f(u),v\rangle
\bigskip $ \\
$\ \ \ \ \ \ \ \ \ \ \ \ \ \ =\langle f_{1}(t),v\rangle -\mu \left(
0,t\right) g_{0}(t)v(0)-\mu \left( 1,t\right) g_{1}(t)v(1),$ $\forall v\in
H^{1},\bigskip $ \\
$u(0)=u_{0}.$%
\end{tabular}%
\right.  \tag{3.44}  \label{c44}
\end{equation}

\textit{Step 4}. \textit{Uniqueness of the solutions}.

First, we shall need the following Lemma.$\bigskip $

\textbf{Lemma 3.5}. \textit{Let} $u$ \textit{be the weak solution of the
following problem}%
\begin{equation}
\left\{
\begin{tabular}{l}
$u_{t}-\frac{\partial }{\partial x}\left[ \mu \left( x,t\right) u_{x}\right]
=\widetilde{f}(x,t),\ 0<x<1,\text{ }0<t<T,$\medskip \\
$u_{x}(0,t)-h_{0}u(0,t)=u_{x}(1,t)+h_{1}u(1,t)=0,$\medskip \\
$u(x,0)=0,$\medskip \\
$u\in L^{2}(0,T;H^{1})\cap L^{\infty }(0,T;L^{2})\cap L^{p}(Q_{T}),$\medskip
\\
$tu\in L^{\infty }(0,T;H^{1}),\text{ \ }tu_{t}\in L^{2}(Q_{T}).$%
\end{tabular}%
\right.  \tag{3.45}  \label{c45}
\end{equation}

\textit{Then}%
\begin{equation}
\begin{tabular}{l}
$\left\Vert u(t)\right\Vert
^{2}+2\int\nolimits_{0}^{t}a(s,u(s),u(s))ds=2\int\nolimits_{0}^{t}\langle
\widetilde{f}(s),u(s)\rangle ds.$%
\end{tabular}
\tag{3.46}  \label{c46}
\end{equation}

The lemma 3.5 is a slight improvement of a lemma used in \cite{1} ( see also
Lions's book \cite{5}).$\blacksquare $

Now, we will prove the uniqueness of the solutions. Assume now that $\left(
H_{7}\right) $ is satisfied.

Let $u_{1}$ and $u_{2}$ be two weak solutions of (\ref{1}) -- (\ref{3}).
Then $u=u_{1}-u_{2}$ is a weak solution of the following problem (\ref{c45})
with the right hand side function replaced by $\widetilde{f}%
(x,t)=-f(u_{1})+f(u_{2}).$ Using Lemma 3.5 we have equality%
\begin{equation}
\begin{tabular}{l}
$\left\Vert u(t)\right\Vert
^{2}+2\int\nolimits_{0}^{t}a(s,u(s),u(s))ds=-2\int\nolimits_{0}^{t}\langle
f(u_{1})-f(u_{2}),u(s)\rangle ds.$%
\end{tabular}
\tag{3.47}  \label{c47}
\end{equation}

Using the monotonicity of $f(y)+\delta y,$ we obtain%
\begin{equation}
\begin{tabular}{l}
$\int\nolimits_{0}^{t}\langle f(u_{1})-f(u_{2}),u(s)\rangle ds\geq -\delta
\int\nolimits_{0}^{t}\left\Vert u(s)\right\Vert ^{2}ds.$%
\end{tabular}
\tag{3.48}  \label{c48}
\end{equation}

It follows from (\ref{c47}), (\ref{c48}) that%
\begin{equation}
\begin{tabular}{l}
$\left\Vert u(t)\right\Vert ^{2}+2a_{0}\int\nolimits_{0}^{t}\left\Vert
u(s)\right\Vert _{H^{1}}^{2}ds\leq 2\delta \int\nolimits_{0}^{t}\left\Vert
u(s)\right\Vert ^{2}ds.$%
\end{tabular}
\tag{3.49}  \label{c49}
\end{equation}

By the Gronwall's Lemma that $u=0.$

Therefore, Theorem 3.1 is proved.$\blacksquare $

\section{\textbf{The boundedness of the solution}}

$\qquad $We now turn to the boundness of the solutions. For this purpose, we
shall make of the following assumptions%
\begin{equation*}
\begin{tabular}{ll}
$\left( H_{1}^{\prime }\right) $ & $h_{0}>0$ and $h_{1}>0,\bigskip $ \\
$\left( H_{2}^{\prime }\right) $ & $u_{0}\in L^{\infty },\bigskip $ \\
$\left( H_{5}^{\prime }\right) $ & $f_{1}\in L^{2}(Q_{T}),\ \ f_{1}(x,t)\leq
0,\text{ \ }a.e.\text{ \ }(x,t)\in Q_{T},\bigskip $ \\
$\left( H_{6}^{\prime }\right) $ & $f\in C^{0}(%
%TCIMACRO{\U{211d} }%
%BeginExpansion
\mathbb{R}
%EndExpansion
)$ satisfies the assumptions $\left( H_{6}\right) ,$ $\left( H_{7}\right) ,$
and$\bigskip $ \\
& $\ \ \ \ \ \ \ \ \ \ \ \ \ \ \ uf(u)\geq 0,\text{ \ }\forall u\in
%TCIMACRO{\U{211d} }%
%BeginExpansion
\mathbb{R}
%EndExpansion
,\text{ \ }\left\vert u\right\vert \geq \left\Vert u_{0}\right\Vert
_{L^{\infty }}.$%
\end{tabular}%
\end{equation*}

We then have the following theorem.

\textbf{Theorem 4.1}. \textit{Let }$(H_{1}^{\prime }),$ $(H_{2}^{\prime }),$
$(H_{3}),$ $(H_{4}),$ $(H_{5}^{\prime }),$ $(H_{6}^{\prime })$ \textit{hold.
Then the unique weak solution of the initial and boundary value problem }(%
\ref{1}) -- (\ref{3}),\textit{\ as given by theorem }3.1\textit{, belongs to}
$L^{\infty }(Q_{T}).$

\textit{Furthermore, we have also}%
\begin{equation}
\begin{tabular}{l}
$\left\Vert u\right\Vert _{L^{\infty }(Q_{T})}\leq \max \left\{ \left\Vert
u_{0}\right\Vert _{L^{\infty }},\text{ }\frac{1}{h_{0}}\left\Vert
g_{0}\right\Vert _{L^{\infty }(0,T)},\frac{1}{h_{1}}\left\Vert
g_{1}\right\Vert _{L^{\infty }(0,T)}\right\} .\blacksquare $%
\end{tabular}
\tag{4.1}  \label{d1}
\end{equation}

\textbf{Remark 4.1}. Assumption $(H_{2}^{\prime })$ is both physically and
mathematically natural in the study of partial differential equation of the
kind of (\ref{1}) -- (\ref{3}), by means of the maximum principle.

\textbf{Proof of Theorem 4.1}. First, let us assume that
\begin{equation}
\begin{tabular}{l}
$u_{0}(x)\leq M,$ $a.e.,\text{\ }x\in \Omega ,$ \ and $\ \max \left\{ \frac{1%
}{h_{0}}\left\Vert g_{0}\right\Vert _{L^{\infty }(0,T)},\frac{1}{h_{1}}%
\left\Vert g_{1}\right\Vert _{L^{\infty }(0,T)}\right\} \leq M.$%
\end{tabular}
\tag{4.2}  \label{d2}
\end{equation}

Then $z=u-M$ satisfies the initial and boundary value%
\begin{equation}
\left\{
\begin{tabular}{l}
$z_{t}-\frac{\partial }{\partial x}\left[ \mu \left( x,t\right) z_{x}\right]
+f(z+M)=f_{1}(x,t),\text{ }0<x<1,\text{ }0<t<T,\bigskip $ \\
$z_{x}(0,t)=h_{0}\left[ z(0,t)+M\right] +g_{0}(t),\text{ \ }-z_{x}(1,t)=h_{1}%
\left[ z(1,t)+M\right] +g_{1}(t),\bigskip $ \\
$z(x,0)=u_{0}(x)-M.$%
\end{tabular}%
\right.  \tag{4.3}  \label{d3}
\end{equation}

Multiplying equation (\ref{d3})$_{1}$ by $v,$ for $v\in H^{1}$ integrating
by parts with respect to variable $x$ and taking into account boundary
condition (\ref{d3})$_{2}$, one has after some rearrangements%
\begin{equation}
\begin{tabular}{l}
$\int\nolimits_{0}^{1}z_{t}vdx+\int\nolimits_{0}^{1}\mu \left( x,t\right)
z_{x}v_{x}dx+\mu \left( 0,t\right) \left[ h_{0}(z(0,t)+M)+g_{0}(t)\text{ }%
\right] v(0)\bigskip $ \\
$\ \ \ \ \ \ \ \ \ \ \ \ \ \ \ \ \ \ \ \ \ \ \ \ \ \ \ \ \ \ \ \ \ \ \ \ \ \
\ +\mu \left( 1,t\right) \left[ h_{1}(z(1,t)+M)+g_{1}(t)\text{ }\right]
v(1)\bigskip $ \\
$\ \ \ \ \ \ \ \ \ \ \ \ \ \ \ \ \ \ \ \ \ \ \ \ \ \ \ \ \ \ \ \ \ \ \ \ \ \
\ +\int\nolimits_{0}^{1}f(z+M)vdx=\int\nolimits_{0}^{1}f_{1}(x,t)vdx,\text{
\ for all \ }v\in H^{1}.$%
\end{tabular}
\tag{4.4}  \label{d4}
\end{equation}

Noticing from assumption $(H_{1}^{\prime })$ we deduce that the solution of
the initial and boundary value problem (\ref{1}) -- (\ref{3}) belongs to $%
L^{2}(0,T;H^{1})\cap L^{\infty }(0,T;L^{2})\cap L^{p}(Q_{T}),$ so that we
are allowed to take $v=z^{+}=\frac{1}{2}(\left\vert z\right\vert +z)$ in (%
\ref{d4}). Thus, it follows that%
\begin{equation}
\begin{tabular}{l}
$\int\nolimits_{0}^{1}z_{t}z^{+}dx+\int\nolimits_{0}^{1}\mu \left(
x,t\right) z_{x}z_{x}^{+}dx+\mu \left( 0,t\right) \left[
h_{0}(z(0,t)+M)+g_{0}(t)\text{ }\right] z^{+}(0,t)\bigskip $ \\
$\ \ \ \ \ \ \ \ \ \ \ \ \ \ \ \ \ \ \ \ \ \ \ \ \ \ \ \ \ \ \ \ \ \ \ \ \ \
\ \ \ +\mu \left( 1,t\right) \left[ h_{1}(z(1,t)+M)+g_{1}(t)\text{ }\right]
z^{+}(1,t)\bigskip $ \\
$\ \ \ \ \ \ \ \ \ \ \ \ \ \ \ \ \ \ \ \ \ \ \ \ \ \ \ \ \ \ \ \ \ \ \ \ \ \
\ \ \
+\int\nolimits_{0}^{1}f(z+M)z^{+}dx=\int%
\nolimits_{0}^{1}f_{1}(x,t)z^{+}dx. $%
\end{tabular}
\tag{4.5}  \label{d5}
\end{equation}

Hence%
\begin{equation}
\begin{tabular}{l}
$\frac{1}{2}\frac{d}{dt}\left\Vert z^{+}(t)\right\Vert
^{2}+a(t,z^{+}(t),z^{+}(t))+\int\nolimits_{0}^{1}f(z^{+}+M)z^{+}dx=\int%
\nolimits_{0}^{1}f_{1}(x,t)z^{+}dx\bigskip $ \\
$\ \ \ \ \ \ \ \ \ \ \ \ \ \ \ \ \ \ \ \ -\mu \left( 0,t\right) \left(
h_{0}M+g_{0}(t)\right) z^{+}(0,t)-\mu \left( 1,t\right) \left(
h_{1}M+g_{1}(t)\right) z^{+}(1,t)\leq 0.$%
\end{tabular}
\tag{4.6}  \label{d6}
\end{equation}%
since%
\begin{equation}
\begin{tabular}{l}
$M\geq \max \{\frac{1}{h_{0}}\left\Vert g_{0}\right\Vert _{L^{\infty }},$ $%
\frac{1}{h_{1}}\left\Vert g_{1}\right\Vert _{L^{\infty }}\}$ and$\bigskip $
\\
$\int\nolimits_{0}^{1}z_{t}z^{+}dx=\int\nolimits_{0,\text{ }%
z>0}^{1}(z^{+})_{t\,}z^{+}dx=\frac{1}{2}\frac{d}{dt}\int\nolimits_{0,\text{
}z>0}^{1}\left\vert z^{+}\right\vert ^{2}dx=\frac{1}{2}\frac{d}{dt}%
\int\nolimits_{0}^{1}\left\vert z^{+}\right\vert ^{2}dx=\frac{1}{2}\frac{d}{%
dt}\left\Vert z^{+}(t)\right\Vert ^{2}.$%
\end{tabular}
\tag{4.7}  \label{d7}
\end{equation}%
and on the domain $z>0$ we have $z^{+}=z$ and $z_{x}=(z^{+})_{x}.$

On the other hand, by the assumption $(H_{2}^{\prime })$ and the inequality (%
\ref{b3}), we obtain%
\begin{equation}
\begin{tabular}{l}
$a(t,z^{+}(t),z^{+}(t))\geq a_{0}\left\Vert z^{+}(t)\right\Vert
_{H^{1}}^{2}. $%
\end{tabular}
\tag{4.8}  \label{d8}
\end{equation}

Using the monotonicity of $f(z)+\delta z$ and $(H_{7})$ we obtain%
\begin{equation}
\begin{tabular}{l}
$\int\nolimits_{0}^{1}f(z^{+}+M)z^{+}dx=\int\nolimits_{0}^{1}\left[
f(z^{+}+M)-f(M)\right] z^{+}dx+\int\nolimits_{0}^{1}f(M)z^{+}dx\bigskip $
\\
\ \ \ \ \ \ \ \ \ \ \ \ \ \ \ \ \ \ \ \ \ \ \ \ \ $\geq -\delta
\int\nolimits_{0}^{1}\left\vert z^{+}\right\vert
^{2}dx+\int\nolimits_{0}^{1}f(M)z^{+}dx\geq -\delta
\int\nolimits_{0}^{1}\left\vert z^{+}\right\vert ^{2}dx=-\delta \left\Vert
z^{+}(t)\right\Vert ^{2}.$%
\end{tabular}
\tag{4.9}  \label{d9}
\end{equation}

Hence, it follows from (\ref{d6}), (\ref{d8}), (\ref{d9}) that%
\begin{equation}
\begin{tabular}{l}
$\frac{1}{2}\frac{d}{dt}\left\Vert z^{+}(t)\right\Vert ^{2}+a_{0}\left\Vert
z^{+}(t)\right\Vert _{H^{1}}^{2}\leq \delta \left\Vert z^{+}(t)\right\Vert
^{2}.$%
\end{tabular}
\tag{4.10}  \label{d10}
\end{equation}

Integrating (\ref{d10}), we get%
\begin{equation}
\begin{tabular}{l}
$\left\Vert z^{+}(t)\right\Vert ^{2}\leq \left\Vert z^{+}(0)\right\Vert
^{2}+2\delta \int\nolimits_{0}^{t}\left\Vert z^{+}(s)\right\Vert ^{2}ds.$%
\end{tabular}
\tag{4.11}  \label{d11}
\end{equation}

Since $z^{+}(0)=(u(x,0)-M)^{+}=(u_{0}(x)-M)^{+}=0,$ hence, using Gronwall's
Lemma, we obtain $\left\Vert z^{+}(t)\right\Vert ^{2}=0.$ Thus $z^{+}=0$ and
$u(x,t)\leq M,$ for $a.e.$ $(x,t)\in Q_{T}.$

The case $-M\leq u_{0}(x),$ $a.e.,\text{\ }x\in \Omega ,$ \ and $M\geq \max
\left\{ \frac{1}{h_{0}}\left\Vert g_{0}\right\Vert _{L^{\infty }(0,T)},\frac{%
1}{h_{1}}\left\Vert g_{1}\right\Vert _{L^{\infty }(0,T)}\right\} $ can be
dealt with, in the same manner as above, by considering $z=u+M$ and $z^{-}=%
\frac{1}{2}(\left\vert z\right\vert -z),$ we also obtain $z^{-}=0$ and hence
$u(x,t)\geq -M,$ for$\ \ a.e.$ $(x,t)\in Q_{T}.$

From all above, one obtains $\left\vert u(x,t)\right\vert \leq M,$ $a.e.$ $%
(x,t)\in Q_{T},$ i.e.,%
\begin{equation}
\begin{tabular}{l}
$\left\Vert u\right\Vert _{L^{\infty }(Q_{T})}\leq M,$%
\end{tabular}
\tag{4.12}  \label{d12}
\end{equation}%
for all $M\geq \max \left\{ \left\Vert u_{0}\right\Vert _{L^{\infty }},\frac{%
1}{h_{0}}\left\Vert g_{0}\right\Vert _{L^{\infty }(0,T)},\frac{1}{h_{1}}%
\left\Vert g_{1}\right\Vert _{L^{\infty }(0,T)}\right\} .$

This implies (\ref{d1}). Theorem 4.1 is proved.$\blacksquare $

\section{\textbf{Asymptotic behavior of the solution as }$t\rightarrow
+\infty .$}

In this part, let $T>0,$ $(H_{1})-(H_{7})$ hold. Then, there exists a unique
solution $u$ of problem (\ref{1}) -- (\ref{3}) such that%
\begin{equation*}
\left\{
\begin{tabular}{l}
$u\in L^{2}(0,T;H^{1})\cap L^{\infty }(0,T;L^{2})\cap L^{p}(Q_{T}),\bigskip $
\\
$tu\in L^{\infty }(0,T;H^{1}),\text{ }\,tu^{\prime }\in L^{2}(Q_{T}).$%
\end{tabular}%
\right.
\end{equation*}

We shall study asymptotic behavior of\ the solution $u(t)$ as $t\rightarrow
+\infty .$

We make the following supplementary assumptions on the functions $\mu \left(
x,t\right) ,$ $f_{1}\left( x,t\right) ,$ $g_{1}(t),$ $g_{2}(t).$%
\begin{equation*}
\begin{tabular}{l}
\begin{tabular}{ll}
$\left( H_{3}^{\prime \prime }\right) $ & $g_{0},$ $g_{1}\in W^{1,1}(%
%TCIMACRO{\U{211d} }%
%BeginExpansion
\mathbb{R}
%EndExpansion
_{+}),\bigskip $%
\end{tabular}
\\
\begin{tabular}{ll}
$\left( H_{4}^{\prime \prime }\right) $ & $\mu \in C^{1}([0,1]\times
%TCIMACRO{\U{211d} }%
%BeginExpansion
\mathbb{R}
%EndExpansion
_{+}),$ $\mu (x,t)\geq \mu _{0}>0,$ $\forall (x,t)\in \lbrack 0,1]\times
%TCIMACRO{\U{211d} }%
%BeginExpansion
\mathbb{R}
%EndExpansion
_{+},\bigskip $%
\end{tabular}
\\
$%
\begin{array}{c}
\left( H_{5}^{\prime \prime }\right)%
\end{array}%
\begin{array}{c}
f_{1}\in L^{\infty }(0,\infty ;L^{2}),%
\end{array}%
\bigskip $ \\
$%
\begin{array}{c}
\left( H_{6}^{\prime \prime }\right)%
\end{array}%
\text{There exist the positive constants }C_{1},$ $\gamma _{1},$ $g_{0\infty
},$ $g_{1\infty }\text{ and the functions}\bigskip $ \\
$\ \ \ \ \ \ \ \ \ \ \ \ \ \
\begin{tabular}{l}
$\mu _{\infty }\in C^{1}([0,1]),\text{ }f_{1\infty }\in L^{2},\text{ such
that}$%
\end{tabular}%
$%
\end{tabular}%
\end{equation*}%
\begin{equation*}
\begin{tabular}{l}
$%
\begin{array}{c}
\begin{tabular}{ll}
$(i)\,$ & $\left\vert g_{0}(t)-g_{0\infty }\right\vert \leq C_{1}e^{-\gamma
_{1}t},\text{ \ }\forall t\geq 0,\bigskip $ \\
$(ii)\,$ & $\left\vert g_{1}(t)-g_{1\infty }\right\vert \leq C_{1}e^{-\gamma
_{1}t},\text{ \ }\forall t\geq 0,\bigskip $ \\
$(iii)\,$ & $\left\Vert \mu \left( t\right) -\mu _{\infty }\right\Vert
_{L^{\infty }}\leq C_{1}e^{-\gamma _{1}t},\text{ \ }\forall t\geq 0,$ $\mu
_{\infty }(x)\geq \mu _{0}>0,\forall x\in \lbrack 0,1],\bigskip $ \\
$(iv)\,$ & $\left\Vert f_{1}(t)-f_{1\infty }\right\Vert \leq C_{1}e^{-\gamma
_{1}t},\text{ \ }\forall t\geq 0.$%
\end{tabular}%
\end{array}%
$%
\end{tabular}%
\end{equation*}

First, we consider the following stationary problem%
\begin{equation}
\left\{
\begin{tabular}{l}
$-\frac{\partial }{\partial x}\left[ \mu _{\infty }\left( x\right) u_{x}%
\right] +f(u)=f_{1\infty }(x),\text{ }0<x<1,\bigskip $ \\
$u_{x}(0)=h_{0}u(0)+g_{0\infty },\text{ }-u_{x}(1)=h_{1}u(1)+g_{1\infty }.$%
\end{tabular}%
\right.  \tag{5.1}  \label{e1}
\end{equation}

The weak solution of problem (\ref{e1}) is obtained from the following
variational problem.

Find $u_{\infty }\in H^{1}$ such that%
\begin{equation}
\begin{tabular}{l}
$a_{\infty }(u_{\infty },v)+\langle f(u_{\infty }),v\rangle =\langle
f_{1\infty },v\rangle -\mu _{\infty }(0)g_{0\infty }v(0)-\mu _{\infty
}(1)g_{1\infty }v(1),$%
\end{tabular}
\tag{5.2}  \label{e2}
\end{equation}%
for all\ $v\in H^{1},$ where%
\begin{equation}
\begin{tabular}{l}
$a_{\infty }(u,v)=\int\nolimits_{0}^{1}\mu _{\infty
}(x)u_{x}(x)v_{x}(x)dx+h_{0}\mu _{\infty }(0)u(0)v(0)+h_{1}\mu _{\infty
}(1)u(1)v(1)\bigskip $ \\
$\ \ \ \ \ \ \ \ \ \ \ \ =\left\langle \mu _{\infty
}u_{x},v_{x}\right\rangle +h_{0}\mu _{\infty }(0)u(0)v(0)+h_{1}\mu _{\infty
}(1)u(1)v(1),$ for all$\mathit{\ }u,v\in H^{1}.$%
\end{tabular}
\tag{5.3}  \label{e3}
\end{equation}

We then have the following theorem.

\textbf{Theorem 5.1}. \textit{Let} $\left( H_{6}\right) ,$ $\left(
H_{3}^{\prime \prime }\right) -\left( H_{6}^{\prime \prime }\right) $
\textit{hold. Then there exists a solution} $u_{\infty }$ \textit{of the
variational problem} (\ref{e2})\textit{\ such that \ }$u_{\infty }\in H^{1}.$

\textit{Furthermore, if }$f$ \textit{satisfies the following condition, in
addition},%
\begin{equation*}
\begin{array}{c}
(H_{7}^{\prime \prime })%
\end{array}%
\text{ \ }%
\begin{tabular}{l}
$f(u)+\delta u\text{ \ \textit{is nondecreasing with respect to variable} }u,%
\text{ \textit{with} }0<\delta <a_{0}.$%
\end{tabular}%
\end{equation*}

\textit{Then the solution is unique}.

\textbf{Proof}. Denote by $\{w_{j}\},$ $\ j=1,2,...$an orthonormal basis in
the separable Hilbert space $H^{1}.$ Put%
\begin{equation}
\begin{tabular}{l}
$y_{m}=$ $\overset{m}{\underset{j=1}{\sum }}d_{mj}w_{j},$%
\end{tabular}
\tag{5.4}  \label{e4}
\end{equation}%
where $d_{mj}$ satisfy the following nonlinear equation system:%
\begin{equation}
\begin{tabular}{l}
$a_{\infty }(y_{m},w_{j})+\langle f(y_{m}),w_{j}\rangle =\langle f_{1\infty
},w_{j}\rangle -\mu _{\infty }(0)g_{0\infty }w_{j}(0)-\mu _{\infty
}(1)g_{1\infty }w_{j}(1),$ $1\leq j\leq m.$%
\end{tabular}
\tag{5.5}  \label{e5}
\end{equation}

By the Brouwer's lemma (see Lions \cite{5}, Lemma 4.3, p. 53), it follows
from the hypotheses $\left( H_{6}\right) ,$ $\left( H_{3}^{\prime \prime
}\right) -\left( H_{6}^{\prime \prime }\right) $ that system (\ref{e4}), (%
\ref{e5}) has a solution $y_{m}.$

Multiplying the $j^{th}$ equation of system (\ref{e5}) by $d_{mj},$ then
summing up with respect to $j,$ we have%
\begin{equation}
\begin{tabular}{l}
$a_{\infty }(y_{m},y_{m})+\langle f(y_{m}),y_{m}\rangle =\langle f_{1\infty
},y_{m}\rangle -\mu _{\infty }(0)g_{0\infty }y_{m}(0)-\mu _{\infty
}(1)g_{1\infty }y_{m}(1).$%
\end{tabular}
\tag{5.6}  \label{e6}
\end{equation}

By using the inequality (\ref{b3}) and by the hypotheses $\left(
H_{6}\right) ,$ $\left( H_{3}^{\prime \prime }\right) -\left( H_{6}^{\prime
\prime }\right) ,$ we obtain%
\begin{equation}
\begin{tabular}{l}
$a_{0}\left\Vert y_{m}\right\Vert _{H^{1}}^{2}+C_{1}\left\Vert
y_{m}\right\Vert _{L^{p}}^{p}\leq C_{1}^{\prime }+\left[ \left\Vert
f_{1\infty }\right\Vert +\sqrt{2}\left( \left\vert \mu _{\infty
}(0)g_{0\infty }\right\vert +\left\vert \mu _{\infty }(1)g_{1\infty
}\right\vert \right) \right] \left\Vert y_{m}\right\Vert _{H^{1}}.$%
\end{tabular}
\tag{5.7}  \label{e7}
\end{equation}

Hence, we deduce from (\ref{e7}) that%
\begin{equation}
\left\{
\begin{tabular}{l}
$\left\Vert y_{m}\right\Vert _{H^{1}}\leq C,\bigskip $ \\
$\left\Vert y_{m}\right\Vert _{L^{p}}\leq C,$%
\end{tabular}%
\right.  \tag{5.8}  \label{e8}
\end{equation}%
$C$ is a constant independent of $m.$

By means of (\ref{e8}) and Lemma 2.1, the sequence $\{y_{m}\}$ has a
subsequence still denoted by $\{y_{m}\}$ such that%
\begin{equation}
\left\{
\begin{tabular}{lll}
$y_{m}\rightarrow u_{\infty }$ & $\text{in}$ & $H^{1}\text{ \ weakly,}$%
\medskip \\
$y_{m}\rightarrow u_{\infty }$ & $\text{in}$ & $L^{2}$ $\text{strongly and \
}a.e.\,\,\text{in \ }\Omega ,$\medskip \\
$y_{m}\rightarrow u_{\infty }$ & $\text{in}$ & $L^{p}\text{ \ weakly.}$%
\end{tabular}%
\right.  \tag{5.9}  \label{e9}
\end{equation}

On the other hand, by (\ref{e9})$_{2}$ and $\left( H_{6}\right) ,$ we have%
\begin{equation}
\begin{tabular}{l}
$f(y_{m})\rightarrow f(u_{\infty })\text{ \ a.e. in \ }\Omega .$%
\end{tabular}
\tag{5.10}  \label{e10}
\end{equation}

We also deduce from the hypothesis $(H_{6})$ and from (\ref{e8})$_{2}$ that%
\begin{equation}
\begin{tabular}{l}
$\int\nolimits_{0}^{1}\left\vert f(y_{m}(x))\right\vert ^{p^{\prime
}}dx\leq 2^{p^{\prime }-1}C_{2}^{p^{\prime
}}[1+\int\nolimits_{0}^{1}\left\vert y_{m}(x)\right\vert ^{p}dx]\leq C,$%
\end{tabular}
\tag{5.11}  \label{e11}
\end{equation}%
where $C$ is a constant independent of $m.$

Applying Lemma 3.4 with $N=1,$ $q=p^{\prime },\ G_{m}=f(y_{m}),$ $%
G=f(u_{\infty }),$ we deduce from (\ref{e10}), (\ref{e11}) that%
\begin{equation}
\begin{tabular}{l}
$f(y_{m})\rightarrow f(u_{\infty })\text{ in }L^{p^{\prime }}\text{ \ weakly}%
.$%
\end{tabular}
\tag{5.12}  \label{e12}
\end{equation}

Passing to the limit in Eq. (\ref{e5}), we find without difficulty from (\ref%
{e9}), (\ref{e12}) that $u_{\infty }$ satisfies the equation%
\begin{equation}
\begin{tabular}{l}
$a_{\infty }(u_{\infty },w_{j})+\langle f(u_{\infty }),w_{j}\rangle =\langle
f_{1\infty },w_{j}\rangle -\mu _{\infty }(0)g_{0\infty }w_{j}(0)-\mu
_{\infty }(1)g_{1\infty }w_{j}(1).$%
\end{tabular}
\tag{5.13}  \label{e13}
\end{equation}

Equation (\ref{e13}) holds for every $j=1,2,...$, i.e., (\ref{e2}) holds.

The solution of the problem (\ref{e2}) is unique; that can be showed using
the same arguments as in the proof of Theorem 3.1.$\blacksquare $

Now we consider asymptotic behavior of\ the solution $u(t)$ as $t\rightarrow
+\infty .$

We then have the following theorem.

\textbf{Theorem 5.2}. \textit{Let} $\left( H_{1}\right) ,$ $\left(
H_{2}\right) ,$ $\left( H_{6}\right) ,$ $\left( H_{3}^{\prime \prime
}\right) -\left( H_{6}^{\prime \prime }\right) ,$ $\left( H_{7}^{\prime
\prime }\right) $ \textit{hold. Then we have}%
\begin{equation}
\begin{tabular}{l}
$\left\Vert u(t)-u_{\infty }\right\Vert ^{2}\leq \left( \left\Vert
u_{0}-u_{\infty }\right\Vert ^{2}+\frac{4C}{\varepsilon \left( \gamma
_{1}-\gamma \right) }\right) e^{-2\gamma t},\text{ }\forall t\geq 0,$%
\end{tabular}
\tag{5.14}  \label{e14}
\end{equation}%
\textit{where}%
\begin{equation*}
\begin{tabular}{l}
$0<\gamma <\min \{\gamma _{1},$ $a_{0}-\delta -4\varepsilon \},$ $\
0<4\varepsilon <a_{0}-\delta ,$%
\end{tabular}%
\end{equation*}%
$C>0$ \textit{is a }a constant independing of $t.$

\textbf{Proof}. Put $Z_{m}(t)=u_{m}(t)-y_{m}.$ Let us subtract (\ref{c5})$%
_{1}$ with (\ref{e5}) to obtain%
\begin{equation}
\left\{
\begin{tabular}{l}
$\langle Z_{m}^{\prime }(t),w_{j}\rangle +a(t;u_{m}(t),w_{j})-a_{\infty
}(y_{m},w_{j})+\langle f(u_{m}(t))-f(y_{m}),w_{j}\rangle \bigskip $ \\
$\ \ \ \ \ \ \ \ \ \ \ \ \ \ \ =\langle f_{1}(t)-f_{1\infty },w_{j}\rangle -%
\left[ \mu \left( 0,t\right) g_{0}(t)-\mu _{\infty }(0)g_{0\infty }\text{ }%
\right] w_{j}(0)\bigskip $ \\
$\ \ \ \ \ \ \ \ \ \ \ \ \ \ \ \ \ \ -\left[ \mu \left( 1,t\right)
g_{1}(t)-\mu _{\infty }(1)g_{1\infty }\text{ }\right] w_{j}(1),$ $1\leq
j\leq m,\bigskip $ \\
$Z_{m}(0)=u_{0m}-y_{m}.$%
\end{tabular}%
\right.  \tag{5.15}  \label{e15}
\end{equation}

By multiplying (\ref{e15})$_{1}$ by $c_{mj}(t)-d_{mj}$ and summing up in $j,$
we obtain%
\begin{equation}
\begin{tabular}{l}
$\frac{1}{2}\frac{d}{dt}\left\Vert Z_{m}(t)\right\Vert
^{2}+a(t;Z_{m}(t),Z_{m}(t))+a(t;y_{m},Z_{m}(t))-a_{\infty
}(y_{m},Z_{m}(t))\bigskip $ \\
$\ \ \ \ \ \ \ \ \ \ \ \ \ \ \ \ \ \ \ +\langle
f(u_{m}(t))-f(y_{m}),Z_{m}(t)\rangle \bigskip $ \\
$\ \ \ \ \ \ \ \ \ \ \ \ \ \ \ =\langle f_{1}(t)-f_{1\infty
},Z_{m}(t)\rangle -\left[ \mu \left( 0,t\right) g_{0}(t)-\mu _{\infty
}(0)g_{0\infty }\text{ }\right] Z_{m}(0,t)\bigskip $ \\
$\ \ \ \ \ \ \ \ \ \ \ \ \ \ \ \ \ \ -\left[ \mu \left( 1,t\right)
g_{1}(t)-\mu _{\infty }(1)g_{1\infty }\text{ }\right] Z_{m}(1,t).$%
\end{tabular}
\tag{5.16}  \label{e16}
\end{equation}

By the assumptions $\left( H_{3}^{\prime \prime }\right) -\left(
H_{6}^{\prime \prime }\right) ,$ $\left( H_{7}^{\prime \prime }\right) ,$
and using the inequalities (\ref{b2}), (\ref{b3}), and with $\varepsilon >0,$
we estimate without difficulty the following terms in (\ref{e16}) as follows%
\begin{equation}
\begin{tabular}{l}
$a(t;Z_{m}(t),Z_{m}(t))\geq a_{0}\left\Vert Z_{m}(t)\right\Vert
_{H^{1}}^{2}; $%
\end{tabular}
\tag{5.17}  \label{e17}
\end{equation}%
\begin{equation}
\begin{tabular}{l}
$\langle f(u_{m}(t))-f(y_{m}),Z_{m}(t)\rangle \geq -\delta \left\Vert
Z_{m}(t)\right\Vert ^{2}\geq -\delta \left\Vert Z_{m}(t)\right\Vert
_{H^{1}}^{2};$%
\end{tabular}
\tag{5.18}  \label{e18}
\end{equation}%
\begin{equation}
\begin{tabular}{l}
$a(t;y_{m},Z_{m}(t))-a_{\infty }(y_{m},Z_{m}(t))=\left\langle \left( \mu
(t)-\mu _{\infty }\right) y_{mx},Z_{mx}(t)\right\rangle \bigskip $ \\
$\ \ \ \ \ \ \ \ \ \ \ \ \ \ \ \ \ \ \ \ \ \ \ \ \ \ \ \ \ \ \ \ \ \ \ \ \ \
\ \ \ \ \ \ \ \ \ \ \ \ +h_{0}\left( \mu \left( 0,t\right) -\mu _{\infty
}(0)\right) y_{m}(0)Z_{m}(0,t)\bigskip $ \\
$\ \ \ \ \ \ \ \ \ \ \ \ \ \ \ \ \ \ \ \ \ \ \ \ \ \ \ \ \ \ \ \ \ \ \ \ \ \
\ \ \ \ \ \ \ \ \ \ \ \ +h_{1}\left( \mu \left( 1,t\right) -\mu _{\infty
}(1)\right) y_{m}(1)Z_{m}(1,t);$%
\end{tabular}
\tag{5.19}  \label{e19}
\end{equation}

Note that $\left\Vert y_{m}\right\Vert _{H^{1}}\leq C,$ we obtain from (\ref%
{e19}) that%
\begin{equation}
\begin{tabular}{l}
$\left\vert a(t;y_{m},Z_{m}(t))-a_{\infty }(y_{m},Z_{m}(t))\right\vert \leq
\left\Vert \mu \left( t\right) -\mu _{\infty }\right\Vert _{L^{\infty
}}\left\Vert y_{mx}\right\Vert \left\Vert Z_{mx}(t)\right\Vert \bigskip $ \\
$\ \ \ \ \ \ \ \ \ \ \ \ \ \ \ \ \ \ \ \ \ \ \ \ \ \ \ \ \ \ \
+2h_{0}\left\Vert \mu \left( t\right) -\mu _{\infty }\right\Vert _{L^{\infty
}}\left\Vert y_{m}\right\Vert _{H^{1}}\left\Vert Z_{m}(t)\right\Vert
_{H^{1}}\bigskip $ \\
$\ \ \ \ \ \ \ \ \ \ \ \ \ \ \ \ \ \ \ \ \ \ \ \ \ \ \ \ \ \ \
+2h_{1}\left\Vert \mu \left( t\right) -\mu _{\infty }\right\Vert _{L^{\infty
}}\left\Vert y_{m}\right\Vert _{H^{1}}\left\Vert Z_{m}(t)\right\Vert
_{H^{1}}\bigskip $ \\
\ \ \ \ \ \ \ \ \ \ \ \ $\leq \left( 1+2h_{0}+2h_{1}\right) C_{1}e^{-\gamma
_{1}t}C\left\Vert Z_{m}(t)\right\Vert _{H^{1}}\leq \varepsilon \left\Vert
Z_{m}(t)\right\Vert _{H^{1}}^{2}+\frac{1}{\varepsilon }Ce^{-2\gamma _{1}t};$%
\end{tabular}
\tag{5.20}  \label{e20}
\end{equation}%
\begin{equation}
\begin{tabular}{l}
$\left\vert \langle f_{1}(t)-f_{1\infty },Z_{m}(t)\rangle \right\vert \leq
\left\Vert f_{1}(t)-f_{1\infty }\right\Vert \left\Vert Z_{m}(t)\right\Vert
\bigskip $ \\
$\ \ \ \ \ \ \ \ \ \ \ \ \ \ \ \ \ \ \ \ \ \ \ \ \ \ \ \leq C_{1}e^{-\gamma
_{1}t}\left\Vert Z_{m}(t)\right\Vert _{H^{1}}\leq \varepsilon \left\Vert
Z_{m}(t)\right\Vert _{H^{1}}^{2}+\frac{1}{\varepsilon }Ce^{-2\gamma _{1}t};$%
\end{tabular}
\tag{5.21}  \label{e21}
\end{equation}%
\begin{equation}
\begin{tabular}{l}
$-\left[ \mu \left( 0,t\right) g_{0}(t)-\mu _{\infty }(0)g_{0\infty }\text{ }%
\right] Z_{m}(0,t)\bigskip $ \\
$\ \ \ \ \ \ \ \ \ \ \ \ \ \ \ \ \ =-\left[ \left( \mu \left( 0,t\right)
-\mu _{\infty }(0)\right) g_{0}(t)+\mu _{\infty }(0)\left(
g_{0}(t)-g_{0\infty }\right) \text{ }\right] Z_{m}(0,t)\bigskip $ \\
$\ \ \ \ \ \ \ \ \ \ \ \ \ \ \ \ \ \leq \sqrt{2}\left\Vert
Z_{m}(t)\right\Vert _{H^{1}}\left[ \left\Vert \mu \left( t\right) -\mu
_{\infty }\right\Vert _{L^{\infty }}\left\Vert g_{0}\right\Vert _{L^{\infty
}(%
%TCIMACRO{\U{211d} }%
%BeginExpansion
\mathbb{R}
%EndExpansion
_{+})}+\mu _{\infty }(0)\left\vert g_{0}(t)-g_{0\infty }\right\vert \right]
\bigskip $ \\
$\ \ \ \ \ \ \ \ \ \ \ \ \ \ \ \ \ \leq \sqrt{2}\left\Vert
Z_{m}(t)\right\Vert _{H^{1}}\left[ \left\Vert g_{0}\right\Vert _{L^{\infty }(%
%TCIMACRO{\U{211d} }%
%BeginExpansion
\mathbb{R}
%EndExpansion
_{+})}+\mu _{\infty }(0)\right] C_{1}e^{-\gamma _{1}t}\leq \varepsilon
\left\Vert Z_{m}(t)\right\Vert _{H^{1}}^{2}+\frac{1}{\varepsilon }%
Ce^{-2\gamma _{1}t}.$%
\end{tabular}
\tag{5.22}  \label{e22}
\end{equation}

Similarly%
\begin{equation}
\begin{tabular}{l}
$-\left[ \mu \left( 1,t\right) g_{1}(t)-\mu _{\infty }(1)g_{1\infty }\text{ }%
\right] Z_{m}(1,t)\leq \varepsilon \left\Vert Z_{m}(t)\right\Vert
_{H^{1}}^{2}+\frac{1}{\varepsilon }Ce^{-2\gamma _{1}t}.$%
\end{tabular}
\tag{5.23}  \label{e23}
\end{equation}

It follows from (\ref{e16}) -- (\ref{e18}), (\ref{e20}) -- (\ref{e23}) and (%
\ref{b3}), that%
\begin{equation}
\begin{tabular}{l}
$\frac{d}{dt}\left\Vert Z_{m}(t)\right\Vert ^{2}+2\left( a_{0}-\delta
-4\varepsilon \right) \left\Vert Z_{m}(t)\right\Vert _{H^{1}}^{2}\leq \frac{8%
}{\varepsilon }Ce^{-2\gamma _{1}t}.$%
\end{tabular}
\tag{5.24}  \label{e24}
\end{equation}

Choose $\varepsilon >0$ and $\gamma >0$ such that $a_{0}-\delta
-4\varepsilon >0$ and $\gamma <\min \{\gamma _{1},$ $a_{0}-\delta
-4\varepsilon \},$ then we have from (\ref{e24}) that%
\begin{equation}
\begin{tabular}{l}
$\frac{d}{dt}\left\Vert Z_{m}(t)\right\Vert ^{2}+2\gamma \left\Vert
Z_{m}(t)\right\Vert ^{2}\leq \frac{8}{\varepsilon }Ce^{-2\gamma _{1}t}.$%
\end{tabular}
\tag{5.25}  \label{e25}
\end{equation}

Hence, we obtain from (\ref{e25}) that%
\begin{equation}
\begin{tabular}{l}
$\left\Vert Z_{m}(t)\right\Vert ^{2}\leq \left[ \left\Vert
Z_{m}(0)\right\Vert ^{2}+\frac{4C}{\varepsilon \left( \gamma _{1}-\gamma
\right) }\right] e^{-2\gamma t}.$%
\end{tabular}
\tag{5.26}  \label{e26}
\end{equation}

Letting $m\rightarrow +\infty $ \ in (\ref{e26}) we obtain%
\begin{equation}
\begin{tabular}{l}
$\left\Vert u(t)-u_{\infty }\right\Vert ^{2}\leq $ $\underset{m\rightarrow
+\infty }{\lim \inf }\left\Vert u_{m}(t)-y_{m}\right\Vert ^{2}\leq \left(
\left\Vert u_{0}-u_{\infty }\right\Vert ^{2}+\frac{4C}{\varepsilon \left(
\gamma _{1}-\gamma \right) }\right) e^{-2\gamma t},\ \text{for all }t\geq 0.$%
\end{tabular}
\tag{5.27}  \label{e27}
\end{equation}

This completes the proof of Theorem 5.2.$\blacksquare $

\section{\textbf{Numerical results}}

$\qquad $First, we present some results of numerical comparison of the
approximated representation of the solution of a nonlinear problem of the
type (\ref{1}) -- (\ref{3}) and the corresponding exact solution of this
problem.

Let the problem%
\begin{equation}
\left\{
\begin{tabular}{l}
$u_{t}-u_{xx}+f(u)=f_{1}(x,t),\text{ }0<x<1,\text{ }t>0,\bigskip $ \\
$u_{x}(0,t)=2u(0,t)+g_{0}(t),\text{ }-u_{x}(1,t)=u(1,t)+g_{1}(t),\bigskip $
\\
$u(x,0)=\widetilde{u}_{0}(x),$%
\end{tabular}%
\right.  \tag{6.1}  \label{f1}
\end{equation}%
where%
\begin{equation}
\left\{
\begin{tabular}{l}
$f_{1}(x,t)=-e^{x}(1+2e^{-t})+(1+e^{-t})^{p-1}e^{(p-1)x},\bigskip $ \\
$f(u)=\left\vert u\right\vert ^{p-2}u,$ $\ p=\frac{5}{2},\bigskip $ \\
$g_{0}(t)=-1-e^{-t},$ $\ g_{1}(t)=-2e(1+e^{-t}),\bigskip $ \\
$\widetilde{u}_{0}(x)=2e^{x}.$%
\end{tabular}%
\right.  \tag{6.2}  \label{f2}
\end{equation}

The exact solution of the problem (\ref{f1}), (\ref{f2}) is $%
u(x,t)=(1+e^{-t})e^{x}.$

To solve numerically the problem (\ref{f1}), (\ref{f2}), we consider the
nonlinear differential system for the unknowns $u_{k}(t)=u(x_{k},t),$ $\
x_{k}=kh,$ $h=1/N.$%
\begin{equation*}
\left\{
\begin{tabular}{l}
$\frac{du_{k}}{dt}(t)=\frac{1}{h^{2}}u_{k-1}-\frac{2}{h^{2}}u_{k}+\frac{1}{%
h^{2}}u_{k+1}-f(u_{k})+f_{1}(x_{k},t),\bigskip $ \\
$u_{0}=\frac{1}{1+2h}\left( u_{1}-hg_{0}(t)\right) ,$ $u_{N}=\frac{1}{1+h}%
\left( u_{N-1}-hg_{1}(t)\right) ,\bigskip $ \\
$u_{k}(0)=\widetilde{u}_{0}(x_{k}),\text{ }k=1,2,...,N-1.$%
\end{tabular}%
\right.
\end{equation*}%
or%
\begin{equation}
\left\{
\begin{tabular}{l}
$\frac{du_{1}}{dt}(t)=\frac{-1}{h^{2}}\left( \frac{1+4h}{1+2h}\right) u_{1}+%
\frac{1}{h^{2}}u_{2}-f(u_{1})-\frac{1}{h(1+2h)}g_{0}(t)+f_{1}(x_{1},t),%
\bigskip $ \\
$\frac{du_{k}}{dt}(t)=\frac{1}{h^{2}}u_{k-1}-\frac{2}{h^{2}}u_{k}+\frac{1}{%
h^{2}}u_{k+1}-f(u_{k})+f_{1}(x_{k},t),\text{ \ }k=\overline{2,N-2},\bigskip $
\\
$\frac{du_{N-1}}{dt}(t)=\frac{1}{h^{2}}u_{N-2}-\frac{1}{h^{2}}\left( \frac{%
1+2h}{1+h}\right) u_{N-1}-f(u_{N-1})-\frac{1}{h(1+h)}%
g_{1}(t)+f_{1}(x_{N-1},t),\bigskip $ \\
$u_{k}(0)=\widetilde{u}_{0}(x_{k}),\text{ \ }k=\overline{1,N-1}.$%
\end{tabular}%
\right.  \tag{6.3}  \label{f3}
\end{equation}

To solve the nonlinear differential (\ref{f3}) at the time $t,$ we use the
following linear recursive scheme generated by the nonlinear term $f(u_{k}):$%
\begin{equation}
\left\{
\begin{tabular}{l}
$\frac{du_{1}^{(n)}}{dt}(t)=\frac{-1}{h^{2}}\left( \frac{1+4h}{1+2h}\right)
u_{1}^{(n)}+\frac{1}{h^{2}}u_{2}^{(n)}-f(u_{1}^{(n-1)})-\frac{1}{h(1+2h)}%
g_{0}(t)+f_{1}(x_{1},t),\bigskip $ \\
$\frac{du_{k}^{(n)}}{dt}(t)=\frac{1}{h^{2}}u_{k-1}^{(n)}-\frac{2}{h^{2}}%
u_{k}^{(n)}+\frac{1}{h^{2}}u_{k+1}^{(n)}-f(u_{k}^{(n-1)})+f_{1}(x_{k},t),\
\text{\ }k=\overline{2,N-2},\bigskip $ \\
$\frac{du_{N-1}^{(n)}}{dt}(t)=\frac{1}{h^{2}}u_{N-2}^{(n)}-\frac{1}{h^{2}}%
\left( \frac{1+2h}{1+h}\right) u_{N-1}^{(n)}-f(u_{N-1}^{(n-1)})-\frac{1}{%
h(1+h)}g_{1}(t)+f_{1}(x_{N-1},t),\bigskip $ \\
$u_{k}^{(n)}(0)=\widetilde{u}_{0}(x_{k}),\text{ \ }k=\overline{1,N-1}.$%
\end{tabular}%
\right.  \tag{6.4}  \label{f4}
\end{equation}

The linear differential system (\ref{f4}) is solved by searching the
associated eigenvalues and eigenfunctions. With a spatial step $h=\frac{1}{5}
$ on the interval $[0,1]$ and for $t\in \lbrack 0,3],$ we have drawn the
corresponding approximate surface solution $(x,t)\longrightarrow u(x,t)$ in
figure 1, obtained by successive re-initializations in $t$ with a time step $%
\Delta t=\frac{1}{50}.$ For comparison in figure 2, we have also drawn the
exact surface solution $(x,t)\longrightarrow u(x,t).$

Note that, the approximate solution $u(x,t)$ decreases exponentially to $%
u_{\infty }(x)$ as $t$ tends to infinity, $u_{\infty }$ being the unique
solution of the corresponding steady state problem%
\begin{equation}
\left\{
\begin{tabular}{l}
$-u_{xx}+\left\vert u\right\vert ^{\frac{1}{2}}u=-e^{x}+e^{\frac{3}{2}x},%
\text{ }0<x<1,\bigskip $ \\
$u_{x}(0)=2u(0)-1,$ $-u_{x}(1)=u(1)-2e.$%
\end{tabular}%
\right.  \tag{6.5}  \label{f5}
\end{equation}%
$\bigskip $

\includegraphics[width=5in]{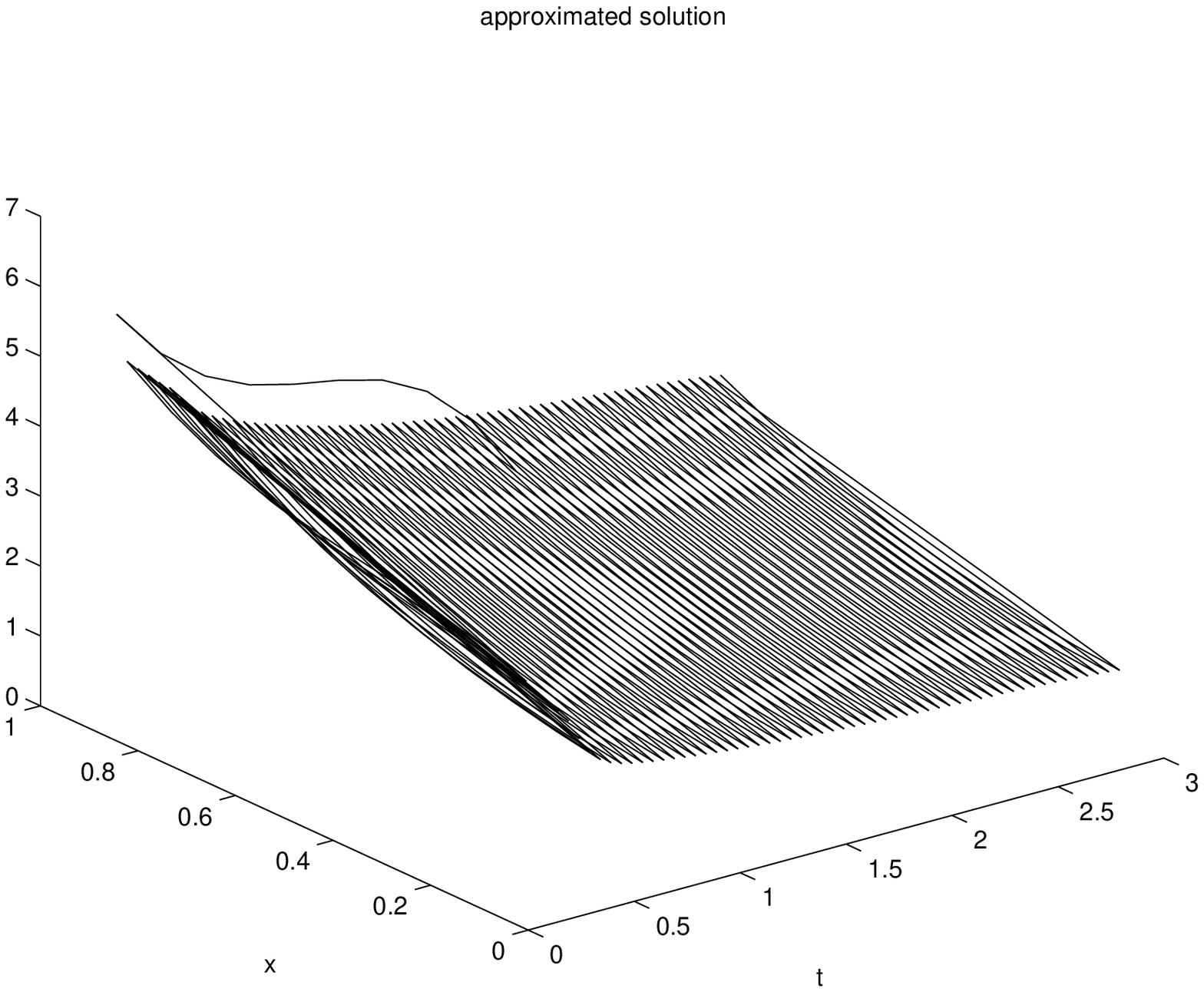}
\begin{center}{Figure 1. Approximated solution} \end{center}

\includegraphics[width=5in]{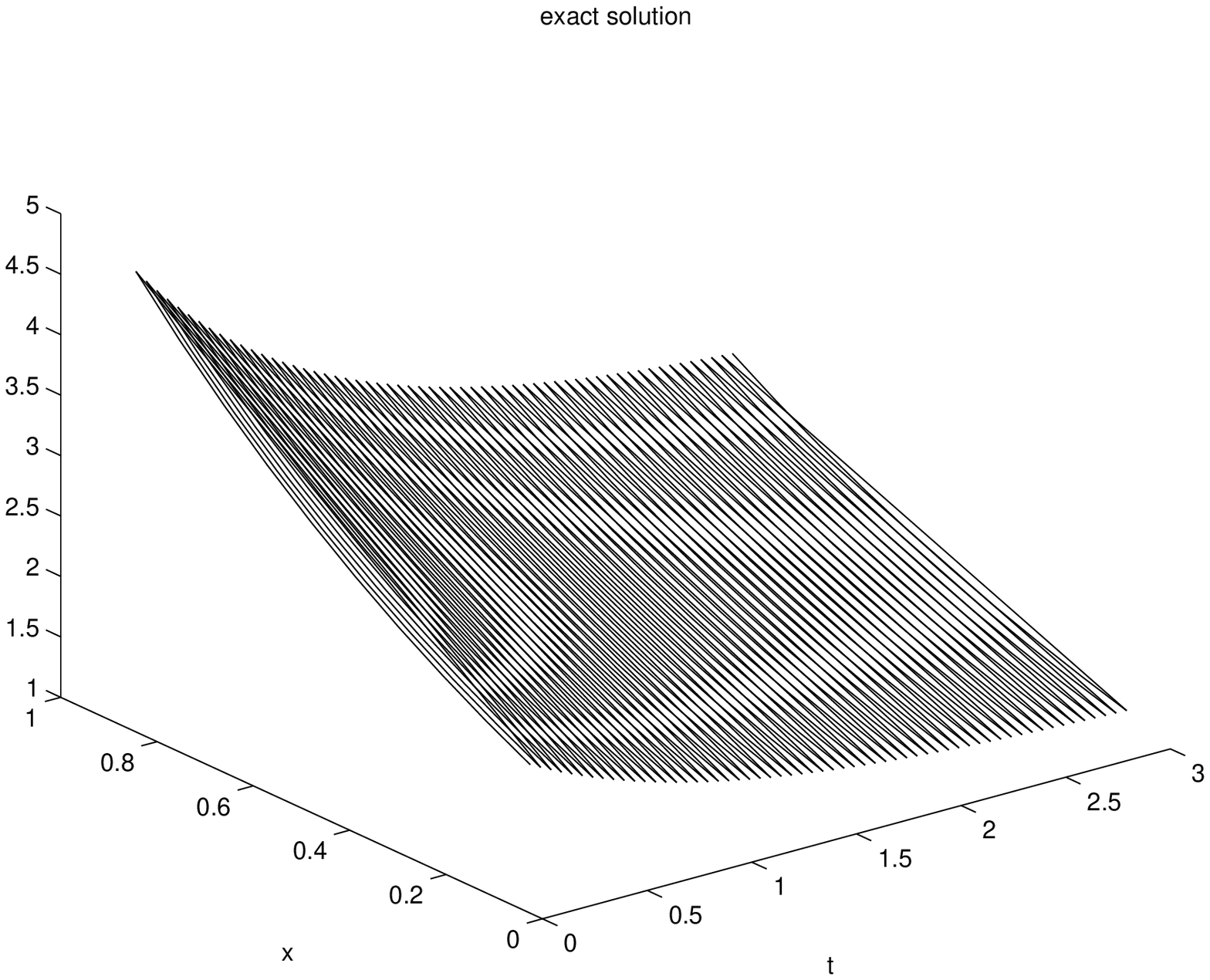}
\begin{center}{Figure 2. Exact solution}\end{center}

\end{document}